 \documentclass[final,3p,times,twocolumn]{elsarticle}

\usepackage{amssymb}
\usepackage{graphicx} 
\usepackage{array} 
\usepackage{verbatim} 
\usepackage{amsmath}
\usepackage{amsfonts}
\usepackage{amssymb}
\usepackage{amsthm}
\newtheorem{thm}{Theorem}[section]
\newtheorem{lem}{Lemma}[section]

\journal{Mathematical Biosciences}

\begin{document}

\begin{frontmatter}



\title{Bifurcation structure of two coupled FHN neurons with delay}

 \author[label1]{Niloofar Farajzadeh Tehrani\corref{cor1}}
 \cortext[cor1]{Corresponding author: Tel.:+982166165626.\\ E-mail address: farajzadeh@mehr.sharif.ir}
 \address[label1]{Department of Mathematical Sciences, Sharif University of Technology, Tehran, Iran.}

\author[label1]{ MohammadReza Razvan}


\begin{abstract}
This paper presents an investigation of the dynamics of two coupled non-identical FitzHugh-Nagumo neurons with delayed synaptic connection. We consider coupling strength and time delay as bifurcation parameters, and try to classify all possible dynamics which is fairly rich.
The neural system exhibits a unique rest point or three ones for the different values of coupling strength by employing the pitchfork bifurcation of non-trivial rest point. The asymptotic
stability and possible Hopf bifurcations of the trivial rest point are studied by analyzing the
corresponding characteristic equation. Homoclinic, fold, 
and pitchfork bifurcations of limit cycles are found. The delay-dependent stability
regions are illustrated in the parameter plane, through
which the double-Hopf bifurcation points can be obtained
from the intersection points of two branches of Hopf
bifurcation. The dynamical behavior of the system
may exhibit one, two, or three different periodic solutions due to pitchfork cycle and torus bifurcations (Neimark-Sacker bifurcation in the Poincare map of a limit cycle), of which detection was impossible without exact and systematic dynamical study. In addition, Hopf, double-Hopf, and torus bifurcations of the non trivial rest points are found. Bifurcation diagrams are obtained numerically or analytically
from the mathematical model and the parameter regions of different behaviors are clarified.
\end{abstract}


\end{frontmatter}


\section{Introduction}
Neurons and their interactions are generally assumed to be the determinant of the brain performance.
The simplest model to display features of neural interactions consists of two coupled neurons or neural systems. 
Starting from such  simple and reduced networks, larger networks can be built
and their features may be studied. In order to study complicated interaction between neurons in large neural networks, the neurons are often put up into highly connected sub-networks or synchronized sub-ensembles. Such neural populations are usually spatially localized \cite{wilson1972excitatory}. In this
way, the model of two mutually coupled neurons may also apply as a framework
of two coupled neural sub-ensembles.\\
  Destexhe et al. \cite{destexhe1994model}, by investigating two interconnected
thalamic spindle oscillatory neurons, showed that how more complex dynamics emanates in ring networks with nearest
neighbors and fully reciprocal connectivity, or in networks in
which every neuron connects to all other nearby neurons. Zhou et al \cite{zhou2006hierarchical} modeled a neural network as a small
sub-network of excitable elements to study synchronization dynamics and the hierarchically clustered
organization of excitable neurons in a network of networks.\\
The FitzHugh-Nagumo (FHN) model with cubic nonlinearity which was derived as a simplified model of the famous Hodgkin-Huxley (HH) model \cite{hodgkin1952quantitative} by FitzHugh \cite{fitzhugh1961impulses} and Nagumo et al. \cite{nagumo1962active}, is a classic oscillator exhibiting variety of nonlinear phenomena in planar autonomous systems. The FHN-like systems are of fundamental importance for describing the qualitative nature of nerve impulse propagation and neural activity.  In fact, this model seems reach enough, and can capture neural excitability of original HH equation \cite{pouryahya2013nonlinear}. In this way, Guckenheimer and Kuehn \cite{guckenheimer2010homoclinic} investigated homoclinic orbits of the FitzHugh–Nagumo equation from the view-point of fast-slow dynamical systems. Hoff et al. \cite{hoff2014numerical} investigated unidirectional and bidirectional  electrical couplings of two identical FHN neurons, involving the effects of the four parameters on bifurcations and synchronization. In modeling studies of trans-membrane potential oscillations, Bachelet et al. \cite{doss2003bursting} explained that bursting oscillations appear quite naturally for two coupled FHN equations. Their coupling model is quite different from other investigations. Different firing patterns such as chaotic firing is found in a pair of identical FHN elements with phase-repulsive coupling \cite{yanagita2005pair}.
 \\
 In addition to study of coupled neurons, the dynamics of a network of FHN neurons is also studied by researchers.
 Oscillatory chaotic behaviors and the mechanism of dynamic change from chaotically spiking to firing death is induced in \cite{ciszak2013experimental}, by considering the coupling of FHN oscillators in a network. Cattani \cite{cattani2014fitzhugh} proposed a network model of instantaneous propagation of signals between neurons in which the neurons need not be physically close to each other. The dynamics of the action potential and the recovery variable within a neural network is described by proposed network. The oscillation patterns in a symmetric network of modified FHN neurons is discussed by Murza \cite{murza2011oscillation}. The author described the building block structure of a three-dimensional lattice shaped as a torus, and identified the symmetry group acting on the electrically coupled differential systems located at the nodes of the lattice.\\
It is known that signal transmission in coupled neurons is not instantaneous in general \cite{swadlow2012axonal}. Hence a time delay can occur in the coupling between neurons or in a self-feedback loop. For example, synaptic communication between neurons depends on the propagation of action potentials across axons. The finite conduction speed and the information processing time in synapses lead to a conduction delay. The speed of signal transmission through unmyelinated axonal fibers is in the order of 1 m/s, which results in existence of  time delays up to 80 ms for propagation through the cortical network \cite{swadlow2012axonal,dhamala2004enhancement}.
 Moreover, the interaction delay between the oscillators can be as large as the oscillation period \cite{izhikevich1998phase}. Therefore, time delay is inevitable in signal transmission for real neurons. On the other hand, time-delayed feedback mechanisms might be intentionally implemented to control neural disturbances, e.g. to suppress alpha rhytm \cite{hadamschek2006brain} and undesired synchrony of firing neurons in Parkinson's disease or epilepsy \cite{schiff1994controlling,rosenblum2004controlling,popovych2005effective}.
Various delayed feedback loops have been proposed as effective and powerful therapy of neurological diseases which are caused by synchronization of neurons \cite{rosenblum2004controlling,popovych2005effective,gassel2007time,gassel2008delay,dahlem2008failure,scholl2009time}. \\
From the modeling viewpoint, delay in dynamical systems is exhibited whenever the system's behavior is dependent at least in part on its history. There are many technological and biological systems which exhibit such behavior. Some examples of delayed systems are coupled lasers  \cite{collera2012bifurcations}, high-speed milling, population dynamics, epidemiology \cite{farajzadeh2013global}, and gene expression, \\
 Krupa and Touboul \cite{krupa2014complex} investigated a self-coupled delayed FHN system to understand the effect of introducing delays in a system whose dynamics, in the absence of delays, is characterized by a canard explosion and relaxation oscillations. They have shown that delays significantly enrich the dynamics, leading to mixed-mode oscillations, bursting, canard and chaos. \\
The research for coupled FHN systems with delay has attracted many authors' attention. Dahlem et al. \cite{dahlem2009dynamics} studied an electrically delay-coupled neural system. They showed that, for sufficiently large delay and coupling strength bi-stability of a fixed point and limit cycle
oscillations occur, even though the single excitable elements display only a stable fixed point.
 Buric and Todorovic \cite{buric2003dynamics}, and also Jia et al. \cite{jia2015dynamic} investigated Hopf bifurcation of coupled FHN neurons with delayed coupling, they concluded that the time delay can induce Hopf bifurcation. Buric et al. \cite{buric2005type}, observed different synchronization states in a coupled FHN system with delay and by  variation of the coupling strength and time delay. Also they showed that the stability and the patterns of exactly synchronous oscillations depend on the type of excitability and type of coupling. In \cite{zhen2010fold} by employing the center manifold theory and normal
form method, the Fold-Hopf bifurcation in a coupled identical FHN neural system with delay is investigated. The synchronization transitions between inphase and antiphase oscillations for delay-coupled FHN system is studied in \cite{song2012inphase}.
 Wang et al. \cite{wang2009bifurcation} reported that time delay can control the occurrence of some bifurcations in two synaptically coupled FHN neurons.  However, only codimension-1 bifurcations are discussed in these papers. Fan and Hong \cite{fan2010hopf}, and also Xu et al. \cite{xu2014dynamics}, considered the stability and Hopf bifurcation of double delay coupled FHN system with different coupling strength. Yao and Tu \cite{yao2014stability} discussed the combined effect of coupling strength and multiple delays on the stability of the rest point, and they obtained stability switches in the coupled FHN neural system with multiple delays. In addition to the study of Hopf bifurcation, Zhen and Xu \cite{zhen2010simple},  and also Rankovic \cite{rankovic2011bifurcations} discussed the effects of the coupling strength and delay on steady state bifurcations in coupled systems with identical neurons. The Bautin bifurcation due to small delay, was analyzed by Buric and Todorovic \cite{buric2005bifurcations}. Nevertheless, in that paper the time-delayed argument of the coupling function was replaced by its Taylor expansion, so that their delay-differential system was approximated by ordinary differential equations. Zhen and Xu \cite{zhen2010bautin} studied Bautin bifurcation of completely synchronous FHN neurons. Li and Jiang \cite{li2011hopf}, investigated Hopf and Bogdanov-Takens bifurcations in a coupled FHN neural system with gap junction for identical neurons. In \cite{panchuk2013synchronization} it was shown that  in a simple motif of two coupled FHN systems with two different delay times and with self-feedback, complex dynamics arise. In particular, resonance effects between the different delay times proved to be crucial.\\
  According to above discussion, FitzHugh-Nagumo model with delayed coupling is interesting and exhibits rich dynamical behavior in comparison with a simplified version of FitzHugh-Nagumo equation \cite{xubifurcations}. Indeed, the dynamics of a pair of excitable systems with time-delayed coupling is quite different from the dynamics with the instantaneous coupling. It is also known that, the behavior of coupled delayed systems differs from that of coupled oscillators with the same coupling \cite{wang2009bifurcation}.\\
In addition to study of coupled FHN neurons, the dynamics of a generalized model consisting of  N  sets of FHN equations has been investigated by Buric and Todorovic \cite{buric2003dynamics}. In their analysis, they applied the terminology found in the description of a small lattice network with a finite time delay between sections for transferring an excitation signal from one element to another. For the case of $N=2$, they studied the properties of stability and the bifurcation of the signals. Kantner et al. \cite{kantner2015delay} showed that how a variety of stable spatio-temporal periodic patterns can be created in 2D-lattices of coupled oscillators with non-homogeneous coupling delays. The effects of heterogeneous coupling delays in complex networks of excitable elements described by the FHN model is studied in \cite{cakan2014heterogeneous}. The effects of discrete as well as of uni- and bi-modal continuous distributions are studied with a focus on regular, small-world, and random networks. The appearance and stability of spatio-temporal periodic patterns in unidirectional rings of coupled FHN systems interacting via excitatory delayed chemical synapses is described in \cite{perlikowski2010periodic}.  A network of connected FHN neurons with delayed coupling and different synaptic strength of self-connection is studied in \cite{lin2014stability}. The stability and oscillation of the solutions is investigated by constructing of Liapunov functional and applying the Chafee’s limit cycle criterion.
Zhen and Xu \cite{zhen2010bautin} and also Lin \cite{lin2013periodic} generalized the FHN models to a three coupled FHN neurons and Feng and Plamondon \cite{feng2014oscillation} to a four coupled FHN neurons with time delay. They derived some criteria to determine the periodic oscillation were provided in the multiple delayed FHN neural system with nonidentical cells.\\
In this paper, we investigate the effect of the coupling strength and time delay on the stability and bifurcations of the system of two coupled FHN neurons. 
Note that neurons in the work of Zhen and Xu \cite{zhen2010simple}, and also Li and Jiang \cite{li2011hopf} are assumed identical which leads to a symmetry for their system of equations, that can not be found in our system. Instead the symmetry in our solutions is due to the inherent symmetry of the FHN neurons, namely our system of equations is equivariant under antipodal map. As a result of the inherent symmetry of the FHN model some pitchfork, pitchfork cycle and Hopf-pitchfork bifurcations are observed in our model.
 We can see that the time delay can either suppress periodic spiking or induce new periodic spiking, depending on the value of the time delay. While in many of the above mentioned papers, only local and codimension-1 bifurcations of FHN systems are studied, the study of global and codimension-2 bifurcations are crucial to understand the mechanism of different dynamical behaviors such as transition between synchronous and anti-phase solutions and simultaneous existence of stable periodic solutions. Actually this kind of bifurcation study is a serious task for the prediction and detection of different dynamical behaviors and special types of solutions. For instance detecting multiple stable periodic solutions which exist simultaneously is rare in literature. In this work we couldn't detect the third stable periodic solution just with simulations. Actually, by the guidance of bifurcation study, we first predicted it's existence, then approved it with simulation.\\
To study the delay effects on neural system in details, in this paper, we have analyzed double-Hopf bifurcations of our system, which is a codimension-two singularity. Various dynamical behaviors are classified in the neighborhood of double-Hopf points. We also show the existence of one, two, or three different modes of spiking as a result of pitchfork cycle and torus (Neimark-Sacker) bifurcations.
We want to emphasize that, by Bendixson-Dulac theorem, it can be seen that the single neuron of the FHN model does not admit periodic solutions for the ranges chosen in our paper, unlike the system in \cite{krupa2014complex} whose dynamics, in the absence of delays, is characterized by a canard explosion and relaxation oscillations. Thus the periodic solutions are due to coupling of the neurons, or time delay. 
Synchronized solutions in many works in the literature are caused by the symmetry of the proposed systems,  like the study of Zhen and Xu \cite{zhen2010simple}, and also Li and Jiang \cite{li2011hopf}. In this paper, we show the existence of both almost synchronized and almost anti-phase solutions for non-identical neurons in different range of parameters. We will show that almost synchronized and almost anti-phase solutions alternatively appear by changing the time delay.\\
 The paper is organized as follows. A brief description of the FHN neuron model is introduced in  Section ~\ref{sec:fhnmodel}. Also  the stability of the trivial rest point is investigated, and the critical values of Hopf and pitchfork bifurcations are driven. In Section ~\ref{sec:small delay}, with the aid of numerical simulations, we have obtained bifurcations for the system with and without delay. Finally, in Section ~\ref{sec:purely} numerical simulations are carried out for determining the bifurcations for a large period of time delay, and delay driven dynamics are analyzed. A conclusion is presented in Section ~\ref{sec:conclusion}.

\section{Model description and basic results}\label{sec:fhnmodel}
In order to consider the effects of delay in the signal transmission between the neurons, we use the paradigmatic FHN model, \cite{fitzhugh1961impulses, nagumo1962active}  with time delay,  investigated by Wang et al. \cite{wang2009bifurcation}.  We here use a coupled nonidentical FHN neural system at which all the parameters are assumed to be non-negative:
 \begin{align}
   \nonumber\dot{v_{1}}(t)&= -v_{1}^{3}(t)+a v_{1}(t)-w_{1}(t)+ c \tanh (v_{2}(t-\tau)), \\
  \nonumber \dot{w_{1}}(t)&=  v_{1}(t)- b_{1} w_{1}(t), \\
    \nonumber\dot{v_{2}}(t)&= -v_{2}^{3}(t)+a v_{2}(t)-w_{2}(t)+c \tanh (v_{1}(t-\tau)), \\
   \dot{w_{2}}(t)&=  v_{2}(t)- b_{2} w_{2}(t), 
\label{1} 
 \end{align}
where $a, b_{1}$ and $b_{2}$ are positive constants, $v_{1}$ represents the
membrane potential, $w_{1}$ is a recovery variable, $ c $ measures the coupling strength,  and $ \tau > 0 $ represents the time delay in signal transmission. We suppose that the function which describes the influence of a neuron on the other one, at time $t $, depends on the state of the neuron at some earlier time $t-\tau$. We consider $c$ and $ \tau $ as bifurcation parameters.
\subsection{Linear stability}
The equations defining the system (\ref{1}) have
$\mathbb{Z}_{2}$-symmetry (invariance under the change of the sign in the state variables). Due to the $\mathbb{Z}_{2}$-symmetry, the origin is always a rest point of the system (\ref{1}), which we call the trivial rest point. The characteristic equation corresponding to linearization of the system (\ref{1}) at the trivial rest point is
 \begin{eqnarray} \label{2}
    P(c,\tau)&=& \lambda^{4}+ A \lambda^{3}+ B \lambda^{2}+ C \lambda + D \\
\nonumber- E (\lambda&+&b_{1})(\lambda+b_{2}) e^{-2\lambda \tau}=0, 
\end{eqnarray}
where $ A=b_{1}+b_{2}-2a $, $ B=b_{1}b_{2}-2a(b_{1}+b_{2})+a^{2}+2 $, $ C= (a^{2}+1)(b_{1}+b_{2})-2ab_{1}b_{2}-2a $, $ D= a^{2}b_{1}b_{2}-a(b_{1}+b_{2})+1$, and $ E= c^{2} $.\\
Due to the presence of the delay, Eq. (\ref{2}) has infinitely many solutions; however, the stability of the fixed points is
determined by a finite number of critical roots with largest real parts.
\subsection{Hopf bifurcation}
We want to obtain some conditions to ensure that the system (\ref{1}) undergoes a single Hopf bifurcation at the trivial rest point $(0,0,0,0)$, when $\tau$ passes through certain critical values.
 Substituting $ \lambda=i\omega$ into Eq. (\ref{2}) and separating the real and imaginary parts, we obtain
 \begin{align} \label{3}
\nonumber  \omega ^{4}&-B \omega^{2}+D-E(-\omega^{2} \cos(2\tau\omega)+b_{1}b_{2}\cos(2\tau\omega)\\
 \nonumber &+(b_{1}+b_{2})\omega sin(2\tau\omega))=0,\\
 \nonumber -A\omega^{3}&+C\omega-E (\omega^{2} \sin(2\tau\omega)-b_{1}b_{2}\sin(2\tau\omega)\\
 & + (b_{1}+b_{2})\omega cos(2\tau\omega))=0.
\end{align} 
Eliminating $\tau$ from Eq. (\ref{3}) gives  
 \begin{eqnarray}  \label{6}
\resizebox{.6\hsize}{!}{$  \omega^{8}+P \omega^{6}+Q \omega^{4}+R \omega^{2}+S=0$},
\end{eqnarray}
where $  P=-2B+A^{2} $, $ Q=B^{2}+2D-2AC-E^{2} $, $ R=-2BD+C^{2}-E^{2}(b_{1}^{2}+b_{2}^{2}) $, and $S= D^{2}-E^{2}(b_{1}b_{2})^{2}$. For $ z=\omega^{2} $, we get
 \begin{eqnarray}\label{4} 
   z^{4}+P z^{3}+Q z^{2}+R z+S=0. 
\end{eqnarray}
Since the form of Eq. (\ref{4}) is identical to that of Eq. (2.4) in the paper of Li and Wei \cite{li2005zeros}, we may apply Lemmas 2.2 and 2.3 analogously. The proofs are omitted.
\begin{lem}\label{lemma1}
If $S<0$, then Eq. (\ref{4}) has at least one positive root.
\end{lem}
Denote $h(z)=z^{4}+P z^{3}+Q z^{2}+R z+S$, then we have $h'(z)= 4 z^{3}+3P z^{2}+2Q z+R$. Set $4 z^{3}+3P z^{2}+2Q z+R=0$. Let $y=z+\frac{3p}{4}$, then the Eq. (\ref{4}) becomes $y^{3}+P_{1}y+Q_{1}=0$, where $P_{1}= \frac{Q}{2}-\frac{3}{16}P^{2}$ and $ Q_{1}=\frac{P^{3}}{32}-\frac{PQ}{8}+\frac{R}{4}$. Define $\Delta=(\frac{Q_{1}}{2})^{2}+(\frac{P_{1}}{3})^{3}$, $\epsilon=\frac{-1+i\sqrt{3}}{2}$,
 \begin{eqnarray}\label{5} 
 \nonumber y_{1}&=&{\textstyle \sqrt[3]{-\frac{Q_{1}}{2}+\sqrt{\Delta}}+ \sqrt[3]{-\frac{Q_{1}}{2}-\sqrt{\Delta}}}, \\
 \nonumber y_{2}&=&{\textstyle\sqrt[3]{-\frac{Q_{1}}{2}+\sqrt{\Delta}}\epsilon+ \sqrt[3]{-\frac{Q_{1}}{2}-\sqrt{\Delta}}}\epsilon^{2}, \\
 \nonumber y_{3}&=&{\textstyle\sqrt[3]{-\frac{Q_{1}}{2}+\sqrt{\Delta}}\epsilon^{2}+ \sqrt[3]{-\frac{Q_{1}}{2}-\sqrt{\Delta}}}\epsilon. 
 \end{eqnarray}
Let $z_{i}=y_{i}-\frac{P}{4}$, ($i$=1,2,3).
\begin{lem}\label{Delta lemma}
Suppose that $ S \geq 0 $, then we have the following results.\\
(i) If $\Delta\geq 0$, then Eq. (\ref{4}) has positive roots if and only if $z_{1}>0$ and $h(z_{1})<0$.\\
(ii) If $\Delta < 0$, then Eq. (\ref{4}) has positive roots if and only if there exists at least one $z^{\ast}\in \{z_{1},z_{2},z_{3}\}$ such that $z^{\ast}>0$ and $h(z^{\ast})\leq 0$.
\end{lem}
Suppose that Eq. (\ref{4}) has positive roots. Without loss of generality, we assume that it has four positive roots, denoted by $z^{\ast}_{k}  (k=1,2,3,4)$. Then Eq. (\ref{6}) has four positive roots, say $\omega_{i}=\sqrt{z^{\ast}_{i}}$, $i=1,2,3,4$. By Eq. (\ref{3}) we have
\small
 \begin{eqnarray}\label{7}
\nonumber  \sin(2\tau \omega)=\frac{(\omega_{k}^{4}-B\omega_{k}^{2}+D)(b_{1}+b_{2})\omega_{k}}{E[(b_{1}+b_{2})^{2}\omega_{k}^{2}+(\omega_{k}^{2}-b_{1}b_{2})^{2}]}\\
\nonumber + \frac{(A\omega_{k}^{3}-C\omega_{k})(-\omega_{k}^{2}+b_{1}b_{2})}{E[(b_{1}+b_{2})^{2}\omega_{k}^{2}+(\omega_{k}^{2}-b_{1}b_{2})^{2}]},\\
\nonumber \cos(2\tau \omega)=\frac{(\omega_{k}^{4}-B\omega_{k}^{2}+D)(b_{1}b_{2}-\omega^{2}_{k}}{E[(b_{1}+b_{2})^{2}\omega_{k}^{2}+(\omega_{k}^{2}-b_{1}b_{2})^{2}]}\\
+\frac{(-A\omega_{k}^{3}+C\omega_{k})((b_{1}+b_{2})\omega_{k})}{E[(b_{1}+b_{2})^{2}\omega_{k}^{2}+(\omega_{k}^{2}-b_{1}b_{2})^{2}]}.
 \end{eqnarray}
\normalsize
Thus, denoting
\small
 \begin{eqnarray} \label{astar}
\nonumber a^{\ast}=\frac{(\omega_{k}^{4}-B\omega_{k}^{2}+D)(b_{1}+b_{2})\omega_{k}}{E[(b_{1}+b_{2})^{2}\omega_{k}^{2}+(\omega_{k}^{2}-b_{1}b_{2})^{2}]}\\
\nonumber + \frac{(A\omega_{k}^{3}-C\omega_{k})(-\omega_{k}^{2}+b_{1}b_{2})}{E[(b_{1}+b_{2})^{2}\omega_{k}^{2}+(\omega_{k}^{2}-b_{1}b_{2})^{2}]},\\
\nonumber b^{\ast}=\frac{(\omega_{k}^{4}-B\omega_{k}^{2}+D)(b_{1}b_{2}-\omega^{2}_{k}}{E[(b_{1}+b_{2})^{2}\omega_{k}^{2}+(\omega_{k}^{2}-b_{1}b_{2})^{2}]}\\
 +\frac{(-A\omega_{k}^{3}+C\omega_{k})((b_{1}+b_{2})\omega_{k})}{E[(b_{1}+b_{2})^{2}\omega_{k}^{2}+(\omega_{k}^{2}-b_{1}b_{2})^{2}]},
 \end{eqnarray}
 \begin{eqnarray}\label{tau}
\nonumber \tau^{(j)}_{k}&=\left\{
\begin{array}{ll}
   \frac{1}{2\omega_{k}}(\arccos b^{\ast}+2j\pi),& a^{\ast}\geqslant 0,\\
    \frac{1}{2\omega_{k}}(2\pi -\arccos b^{\ast}+2j\pi),& a^{\ast}< 0,
\end{array} \right.\\
 &(k=1,2,3,4, j=0,1,2,...)
  \end{eqnarray}
\normalsize
then $\pm i\omega_{k} $ is a pair of purely imaginary roots of (\ref{2}) with $\tau=\tau^{(j)}_{k}$.\\
Applying Lemmas \ref{lemma1}, \ref{Delta lemma}, the Routh-Hurwitz (R-H) criterion, and Ruan and Wei's result \cite{ruan2003zeros} according to Rouche’s theorem, we have the following results.
\begin{lem}\label{stability}
Assume that $A>0$, $A(B-E)>C-E(b_{1}+b_{2})$, $D>Eb_{1}b_{2}$, and $[C-E(b_{1}+b_{2})][A(B-E)-C+E(b_{1}+b_{2})]>A^2(D-Eb_{1}b_{2})$ are satisfied, ((R-H) hypothesis).\\
(i) If one of the followings holds: (a) $S<0$; (b) $S\geq 0, D\geq 0, z_{1}>0$, and $h(z_{1})\leq 0$; (c) $S\geq0$, $D<0$, and there exists $z^{*}\in \{z_{1},z_{2},z_{3}\}$ such that $z^{*}>0$ and $h(z^{*})\leq 0$, then all roots of (\ref{2}) have negative real parts when $\tau \in [0,\tau_{0})$, such that $\tau_{0}=min \{\tau^{(0)}_{1},\tau^{(0)}_{2},\tau^{(0)}_{3},\tau^{(0)}_{4}\}$.\\
(ii) If the conditions (a)-(c) of (i) are not satisfied, then all roots of (\ref{2}) have negative real parts for all $\tau \geq 0$.
\end{lem}
Motivated by Lemma 2.3 of the work of Ruan and Wei \cite{ruan2003zeros}, Lemma 2.4 of the work of Li and Wei \cite{li2005zeros}, Lemma 2.5 of the paper of Hu and Huang \cite{hu2009stability}, and  Theorem 2.1 of Fan and Hong \cite{fan2010hopf}, we obtain following conclusions.
\begin{lem}\label{realpart}
Suppose $h'(z_{0})\neq 0$. If $\tau=\tau_{0}$, then $\pm i\omega_{0}$ is a pair of simple purely roots of Eq. (\ref{2}). In addition, $\frac{(d\Re\lambda(\tau))}{d\tau}\vert_{\tau=\tau^{(j)}_{k}}\neq 0$, and the sign of $\frac{(d\Re\lambda(\tau))}{d\tau}\vert_{\tau=\tau^{(j)}_{k}}$ is consistent with that of $h'(z^*_{k})$.
\end{lem}
Applying Lemmas \ref{stability}-\ref{realpart}, we obtain the following theorem immediately.
\begin{thm}
Suppose hypothesis (R-H) of Lemma \ref{stability} hold.\\
(i) If non of the conditions (a) $S<0$; (b) $S\geq 0, D\geq 0, z_{1}>0$, and $h(z_{1})\leq 0$; (c) $S\geq0$, $D<0$, and there exists a $z^{*}\in \{z_{1},z_{2},z_{3}\}$ such that $z^{*}>0$ and $h(z^{*})\leq 0$ is satisfied, then the zero solution  of (\ref{1}) is asymptotically stable for all $\tau\geq0$.\\
(ii) If one of the conditions (a), (b), or (c) of (i) is satisfied, then the zero solution of system (\ref{1}), for $\tau \in [0,\tau_{0})$, is asymptotically stable ( $\tau_{0}$ is the parameter defined by the Lemma \ref{stability}).\\
(iii) If one of the conditions (a), (b), and (c) of (i) is satisfied, and  $h'(z^*_{k})\neq 0$, then for $\tau=\tau^{(i)}_{k}, (i=1,2,3,...)$, the system (\ref{1}) undergoes a Hopf bifurcation at $(0,0,0,0)$.
\end{thm}
\subsection{Pitchfork bifurcation}
Now we want to study the possible steady state bifurcations of the trivial rest point.
Motivated by Lemma 2.1 of the work of Li and Jiang  \cite{li2011hopf}, we have the following conclusions on the eigenvalues of the system (\ref{1}).
 \begin{lem}
 Eq. (\ref{2}) has a zero eigenvalue if and only if $ c^{2}=\frac{a^{2}b_{1}b_{2}+1-a(b_{1}+b_{2})}{b_{1}b_{2}}=\frac{D}{b_{1}b_{2}} $ and $ \tau\neq \frac{2ab_{1}b_{2}+2a-(a^{2}+1-c^{2})(b_{1}+b_{2})}{2c^{2}b_{1}b_{2}}=\frac{c^{2}-C}{2c^{2}b_{1}b_{2}} $.\\
 $\lambda=0$ is a double root of (\ref{2}), if and only if $ c^{2}=\frac{a^{2}b_{1}b_{2}+1-a((b_{1}+b_{2}))}{b_{1}b_{2}}=\frac{D}{b_{1}b_{2}} $, $ \tau=\frac{2ab_{1}b_{2}+2a-(a^{2}+1-c^{2})(b_{1}+b_{2})}{2c^{2}b_{1}b_{2}}=\frac{c^{2}-C}{2c^{2}b_{1}b_{2}} $, and
 \begin{align*} 
  \nonumber \resizebox{.95\hsize}{!}{$\frac{\partial^{2} P(c,\tau) }{ \partial^{2} \lambda}\vert_{\lambda=0}=  b_{1}b_{2}(1-2 c^{2}\tau^{2})+a^{2}-c^{2}+2+(b_{1}+b_{2})(2\tau c^{2}-2a)$}.
 \end{align*}
\label{pitchfork-hopf}
 \end{lem}
 We can see that our model (\ref{1}) have $\mathbb{Z}_{2}$-symmetry, and in the presence of this symmetry the pitchfork bifurcation becomes generic \citep{ruelle2014elements}. Therefore, possible zero eigenvalue of the Eq. (\ref{2}) is correspondent to pitchfork bifurcation, and saddle-node bifurcation is impossible. Moreover, possible double zero eigenvalue of the Eq. (\ref{2}) is correspondent to Hopf-pitchfork  bifurcation.\\
 In the rest of the paper we will explain coupling and delay driven dynamics of the system (\ref{1}). Bifurcation diagrams are obtained numerically or analytically from the mathematical model and the parameter regions of different behavior are clarified. Numerical simulations using the bifurcation analysis software DDE-Biftool \cite{engelborghs2002numerical} are carried out to illustrate the validity of the main results. 
\section{The dynamics for small time delays}\label{sec:small delay}
We want to consider possible rest points and limit cycles of the system (\ref{1}), and their bifurcations for the variable parameters $c$ and small $\tau$'s. As stated in the previous section the origin is always a rest point of the system, without any constraint on the parameters. In the paper we choose the value of the parameters as $a=0.55$, $b_{1}=1.128$, and $b_{2}=0.58$. This choice of values makes the trivial rest point stable for $c=0$. Hence, by this choice of the parameters the neurons are at rest, without coupling. First, we want to study the coupling strength-driven dynamics of the system. Therefore, as the next step, let us briefly consider the system (\ref{1}), for  $\tau = 0$.
\subsection{Instantaneous coupling}\label{sec:instantaneous}
We know that $(0,0,0,0)$ is always a rest point of the system (\ref{1}), and according to the $\mathbb{Z}_{2}$-symmetry of the system, if $(v_{1}^{*},w_{1}^{*},v_{2}^{*},w_{2}^{*})$ is a rest point of the system (\ref{1}), so is $-(v_{1}^{*},w_{1}^{*},v_{2}^{*},w_{2}^{*})$. Therefore, pitchfork bifurcation is generic for our model, and the double zero eigenvalue of the Eq. (\ref{2}) is correspondent to pitchfork bifurcation. For system (\ref{1}) we can see that for fixed parameters $a=0.55$, $b_{1}=1.128$, $b_{2}=0.58$, and for the parameter $c=0.6285$, the pitchfork bifurcation occurs and a pair of new rest points appear.
\begin{figure}
\includegraphics[width=0.5\textwidth]{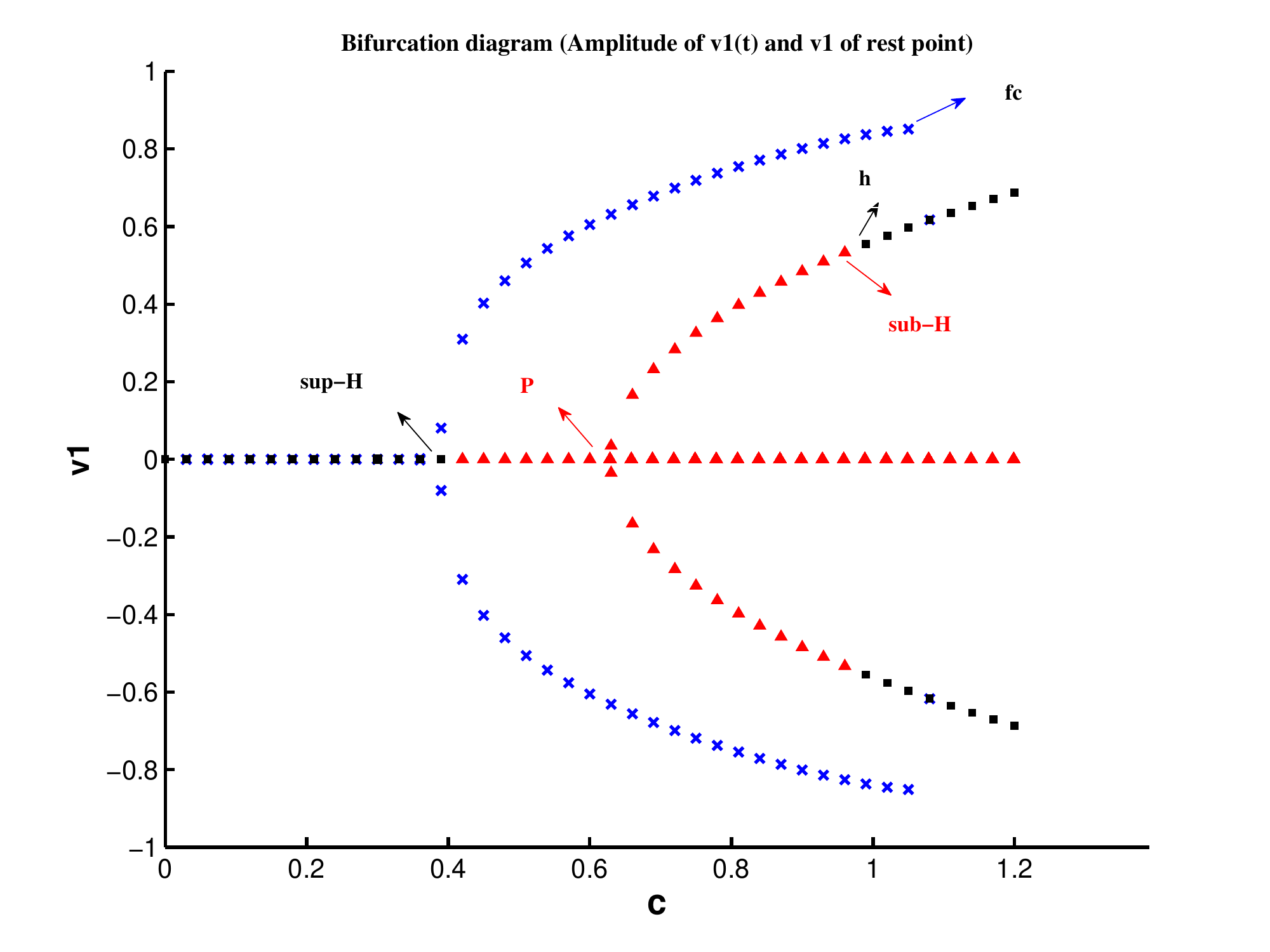} 
\caption{\footnotesize The squared (black) line shows the stable rest points of coupled neurons without delay. For $c<0.3974$ the origin is the only rest point. The cross (blue) line shows the amplitude of the periodic cycle which appears after super-critical Hopf bifurcation (sup-H) of the trivial rest point. After Hopf bifurcation the trivial rest point becomes unstable, which is shown by triangle marker. When pitchfork bifurcation (P) occurs at $c=0.6285$, a pair of nontrivial rest points appear which are unstable and are shown by triangle marker. The nontrivial rest points become stable due to sub-critical Hopf bifurcation (sub-Hopf) at $c=0.9751$, and two unstable limit cycles appear. Two unstable limit cycles collide at $c=1.0545$, and homoclinic bifurcation occurs (h). The stable and unstable limit cycles annihilate each other and a fold bifurcation of limit cycles (fc) occurs at $c=1.0721$.}
\label{pitchfork}
\end{figure}  
 To fully study the bifurcations of the trivial rest point we will increase the parameter $c$ continuously, Fig. \ref{pitchfork}. If we consider two coupled neurons without delay, for small amount of the parameter $c$ the trivial rest point is stable. By increasing the parameter $c$, it is seen that for $c=0.3974$ a super-critical Hopf bifurcation occurs, Fig. \ref{pitchfork}(sup-H),  and the neuron's stability switches from resting state to periodic spiking. The amplitude of the stable limit cycle grows continuously with $c$, in a small interval of the values of $c$. By increasing the parameter $c$, there is a sub-critical pitchfork bifurcation at the parameter $c=0.6285$, Fig. \ref{pitchfork}(P),  and a pair of antipodal rest points appear. The stability of the origin before and after the pitchfork bifurcation and the two new rest points depends on the occurrence of Hopf bifurcation for the trivial rest point. Since the super-critical Hopf bifurcation occurs before pitchfork bifurcation, all the three rest points of the system are unstable after pitchfork bifurcation, and the system is mono-stable, i.e. the only stable state of the system is the limit cycle that stated above, see Fig.  \ref{pitchfork}.  
By increasing the parameter $c$, for $c=0.9751$,  for two non-trivial rest points a sub-critical Hopf bifurcation occurs, and two small unstable limit cycles appear around non trivial rest points, Fig. \ref{pitchfork}(sub-H). Hence, two nontrivial rest points become stable and the system becomes dichotomous; two rest points and a limit cycle. By increasing the parameter $c$, two unstable limit cycles grow continuously until they collide at $c=1.0545$, Fig. \ref{pitchfork}(h). Therefore, as a global bifurcation, an eight-figure saddle homoclinic  bifurcation occurs. 
 At the bifurcation point the period of the unstable limit cycles have grown to infinity, see Fig. \ref{infinite period}(Left).  After this bifurcation a big unstable limit cycle around the rest points appears. We should emphasize that after this bifurcation the system remains dichotomous and the behavior of the neurons doesn't change. It's not hard to check that the saddle quantity is positive and according to the theorems on homoclinic orbits developed by Shilnikov \cite{kuznetsov2013elements}, the bifurcation results in the birth of a unique unstable limit cycle from homoclinic orbit.  The unstable eight-figure saddle homoclinic orbit is shown in Fig. \ref{infinite period}(Right).
 \begin{figure}
\begin{tabular}{l}
\includegraphics[width=4.5cm]{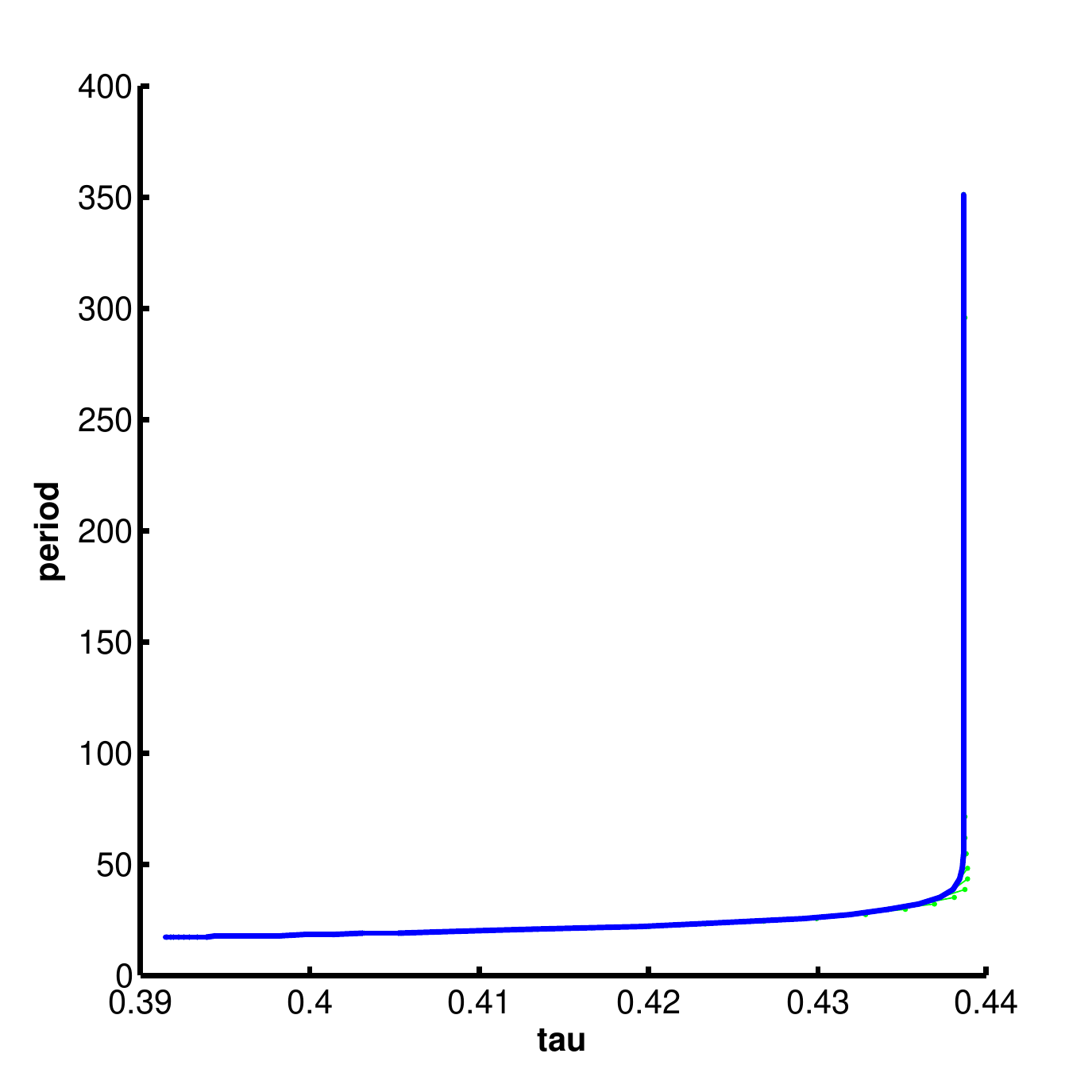}
\includegraphics[width=4.5cm]{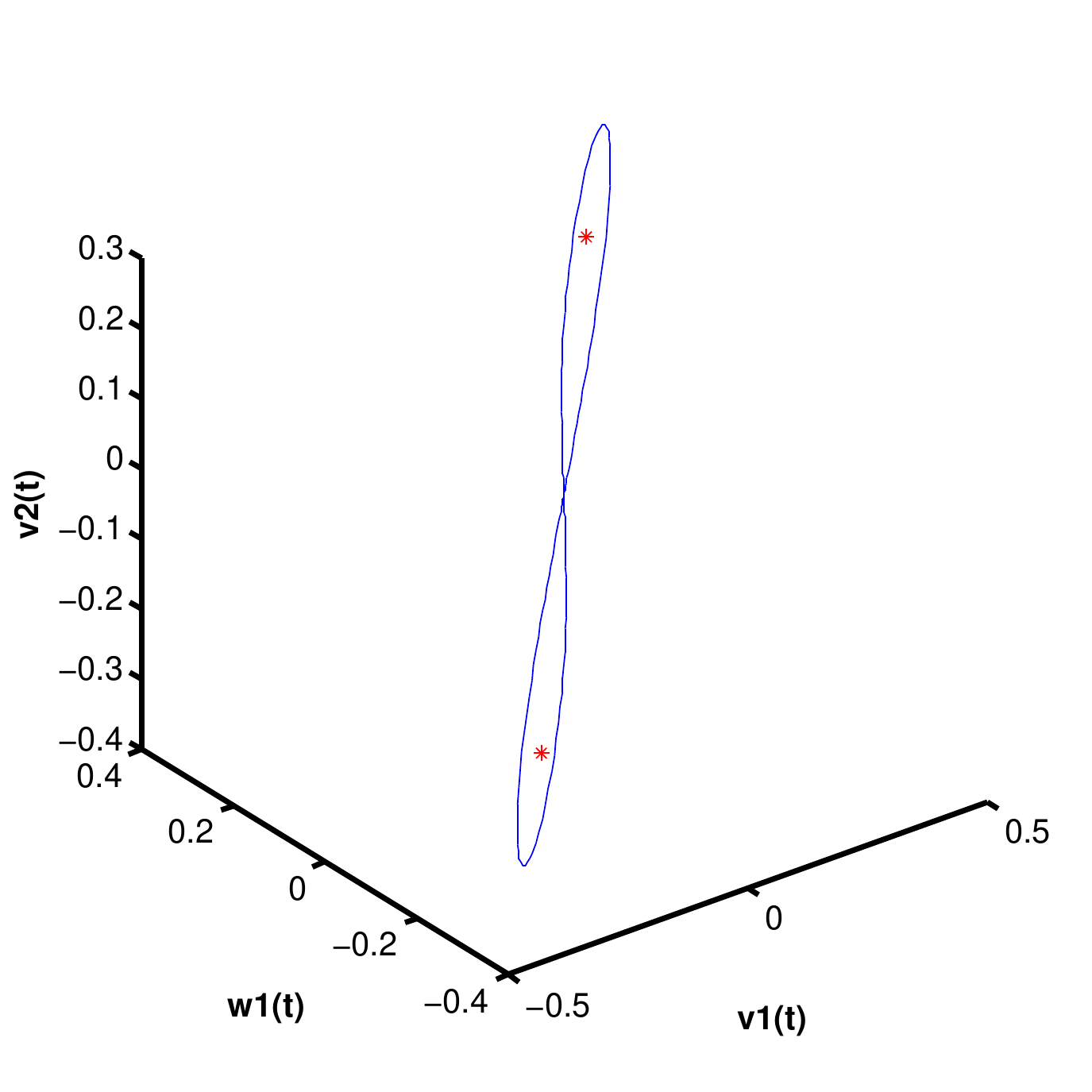}
  \end{tabular}
   \caption{\footnotesize (Left) The period of the unstable limit cycle which appears through the sub-critical Hopf bifurcation. At the bifurcation point $c=1.0545$, the period of the unstable limit cycle has grown to infinity. (Right) The unstable eight-figure saddle homoclinic orbit. } 
 \label{infinite period}
\end{figure}
 By a small change of the parameter $c$, at $c\simeq 1.0721$ a fold of limit cycles bifurcation occurs,  point (fc) in Fig. \ref{pitchfork}. Therefore, the stable limit cycle emanating from the Hopf bifurcation of the trivial rest point and the new born unstable limit cycle approach and annihilate each other, see Fig. \ref{fig3}. Thus the only stable states of the system are non trivial rest points and the system becomes quiescent. 
\begin{figure}
\includegraphics[width=0.4\textwidth]{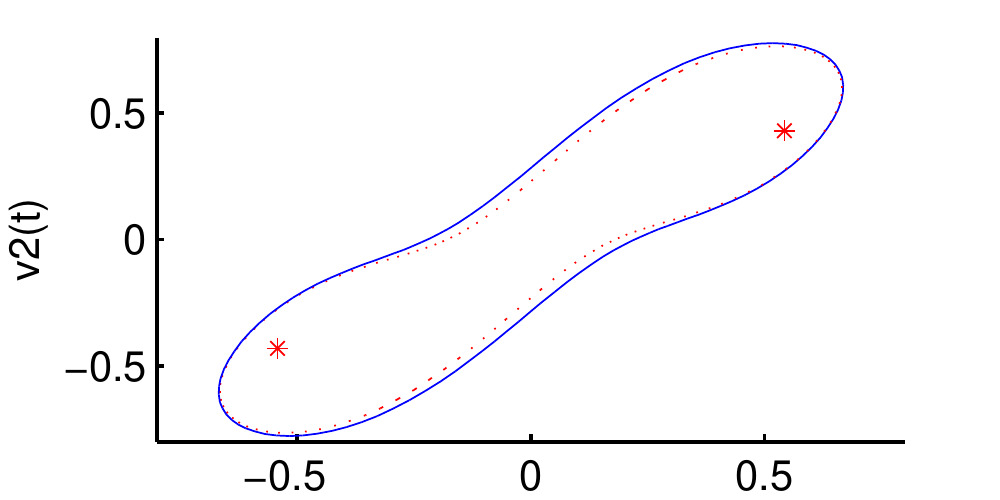}
\caption{\footnotesize Stable (blue solid line) and unstable (red dashed) limit cycles near the fold of limit cycles bifurcation.}
\label{fig3}
\end{figure}
\subsection{Delayed coupling}
Now we come back to the system (\ref{1}) with small delays. First we want to consider the bifurcations of the trivial rest point. From now on, the bifurcation parameters are coupling strength $c$ and time delay  $\tau$.
\subsubsection{Bifurcations of the trivial rest point}
\begin{figure}
\centering
\includegraphics[width=0.5\textwidth]{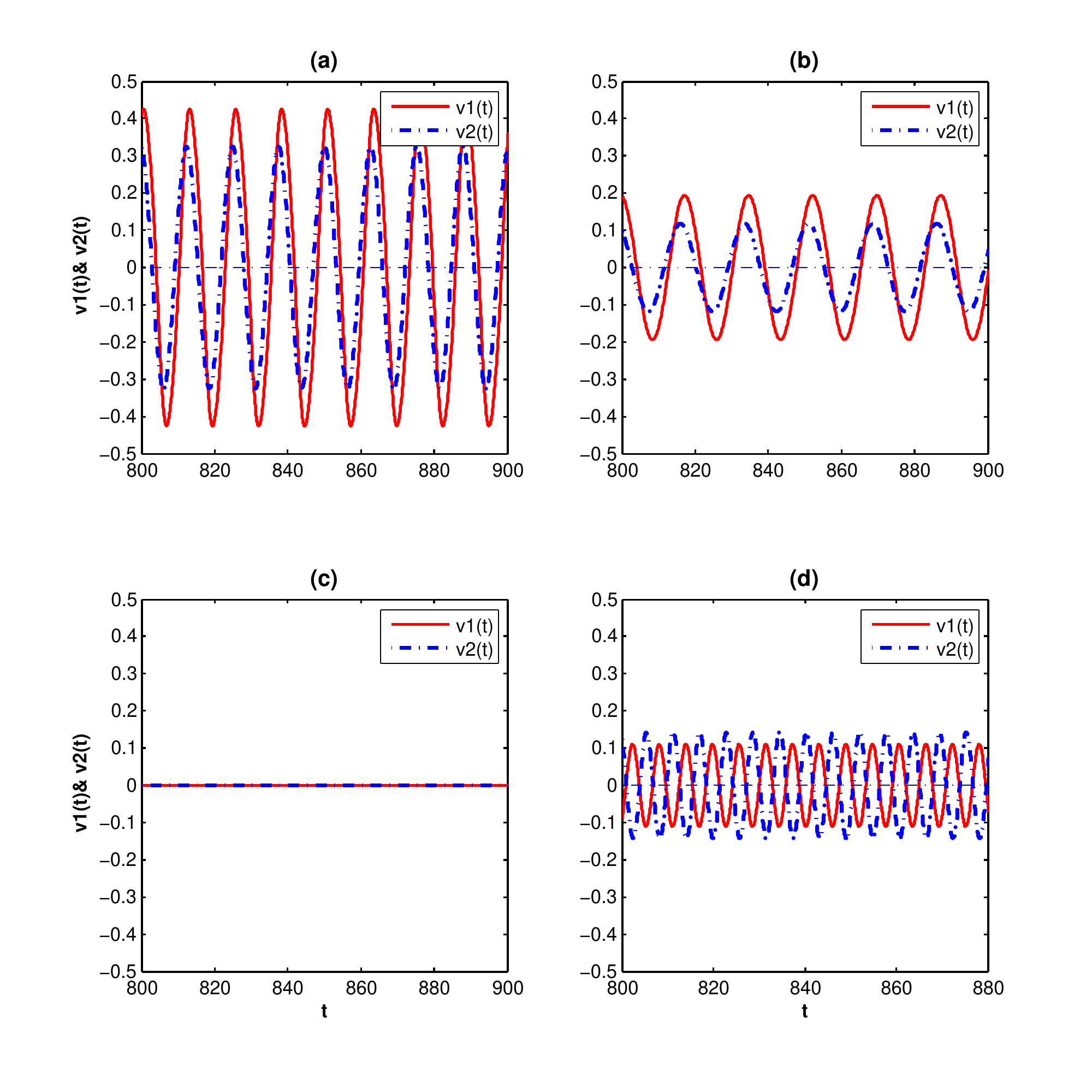}
\caption{\footnotesize The impact of the time delay on the dynamics of the neurons for $c=0.5$.  (a) $\tau=0.1$, the neurons spike with a large amplitude. (b) $\tau=0.25$, the neurons spike with a small amplitude and are near the super-critical Hopf bifurcation. (c) $\tau=0.4$, the neurons are quiescence. (d) $\tau=1.2$, the neurons are spiking.}
\label{impactoftau}
\end{figure}

The type of bifurcations of the rest points determines the most important neuro-computational features of the neuron. For the system (\ref{1}) neurons are excitable because the trivial resting state is near a Hopf bifurcation, i.e., near a transition from quiescence to spiking. Hence, by increasing  $c$ and $\tau$, a branch of Hopf bifurcation emanates from $c=0.3974$ and $\tau=0$ in $(c,\tau)$ plane, and continues until $c=0.6285$ and $\tau=0.5219$ (first br of Hopf in Fig. \ref{fig8}), in which the Hopf
curve meets the curve corresponding to the pitchfork bifurcation of the trivial rest point. In this situation the second conditions of Lemma \ref{pitchfork-hopf} are satisfied and a Hopf-pitchfork bifurcation occurs, see Fig. \ref{fig8}. Since the Hopf bifurcation occurs before pitchfork bifurcation, the two nontrivial rest points are both unstable.
 To illustrate what we have stated above, we can fix the parameter $c$ and study the dynamic changes according to $\tau$, see Fig. \ref{impactoftau}.
As it is shown in Fig. \ref{impactoftau}(a), for $c=0.5$ and $\tau=0.1$  the neurons oscillate almost synchronized, in  which the
phase shift between oscillation of two neurons is close (but not equal) to zero. By increasing the parameter $\tau$ the amplitude of the periodic solutions decrease until the parameter $\tau$ reaches the first branch of Hopf bifurcation. Due to the super-critical Hopf bifurcation the stable limit cycle disappears and the trivial rest point becomes stable, Fig. \ref{impactoftau}(c). By further increasing the parameter $\tau$, the neurons start to spike in an almost anti-phase manner, in which the
phase shift between oscillation of two neurons  is close (but not equal) to $\pi$, Fig. \ref{impactoftau}(d). We will see in the next sections that, by increasing the range of parameter $\tau$ other branches of Hopf bifurcation appear which are purely delay-driven and have no trace in system with instantaneous coupling.
We should notice that, even for small delays, the parameter $\tau$ can either suppress periodic spiking or induce new periodic spiking, depending on the value of the time delay, see Fig. \ref{impactoftau}. Thus, treatment of neural systems with changing the time delay can be delicate and challenging.
\subsubsection{Bifurcations of the non-trivial rest points}
\begin{figure}
\includegraphics[width=0.4\textwidth]{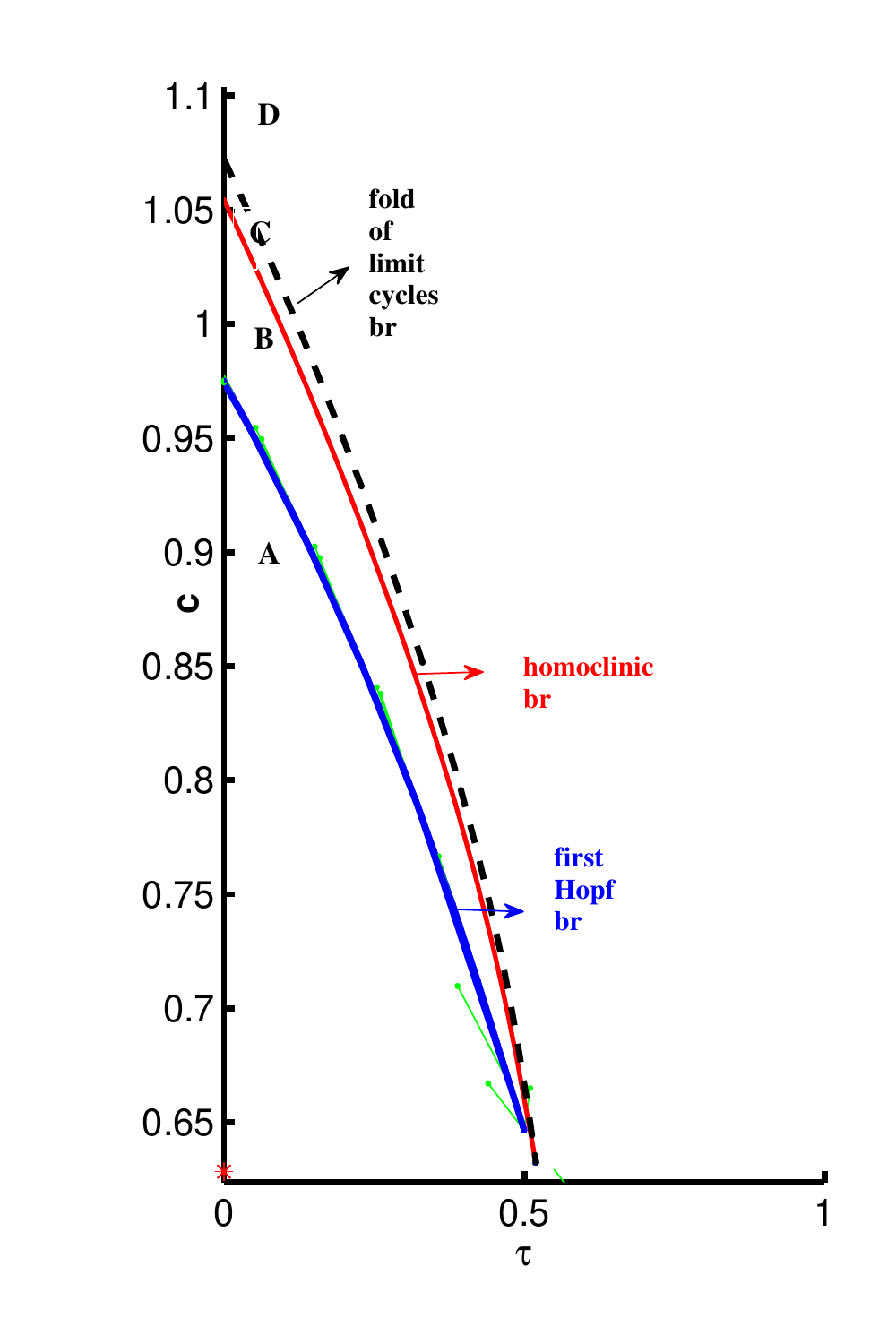}
\caption{\footnotesize Bifurcation diagram. Blue line is a branch of Hopf bifurcation for nontrivial rest points for small delays. The red line is eight-figure saddle homoclinic branch which emanates at $c=1.0545$.  The dashed black line is the branch of fold of limit cycles, which emanates at $c=1.0721$. The (red) stared point, shows the pitchfork bifurcation when $\tau=0$. The solid green line is the branch of pitchfork bifurcation. }
\label{fig6}
\end{figure}

\begin{figure}
\includegraphics[width=0.5\textwidth]{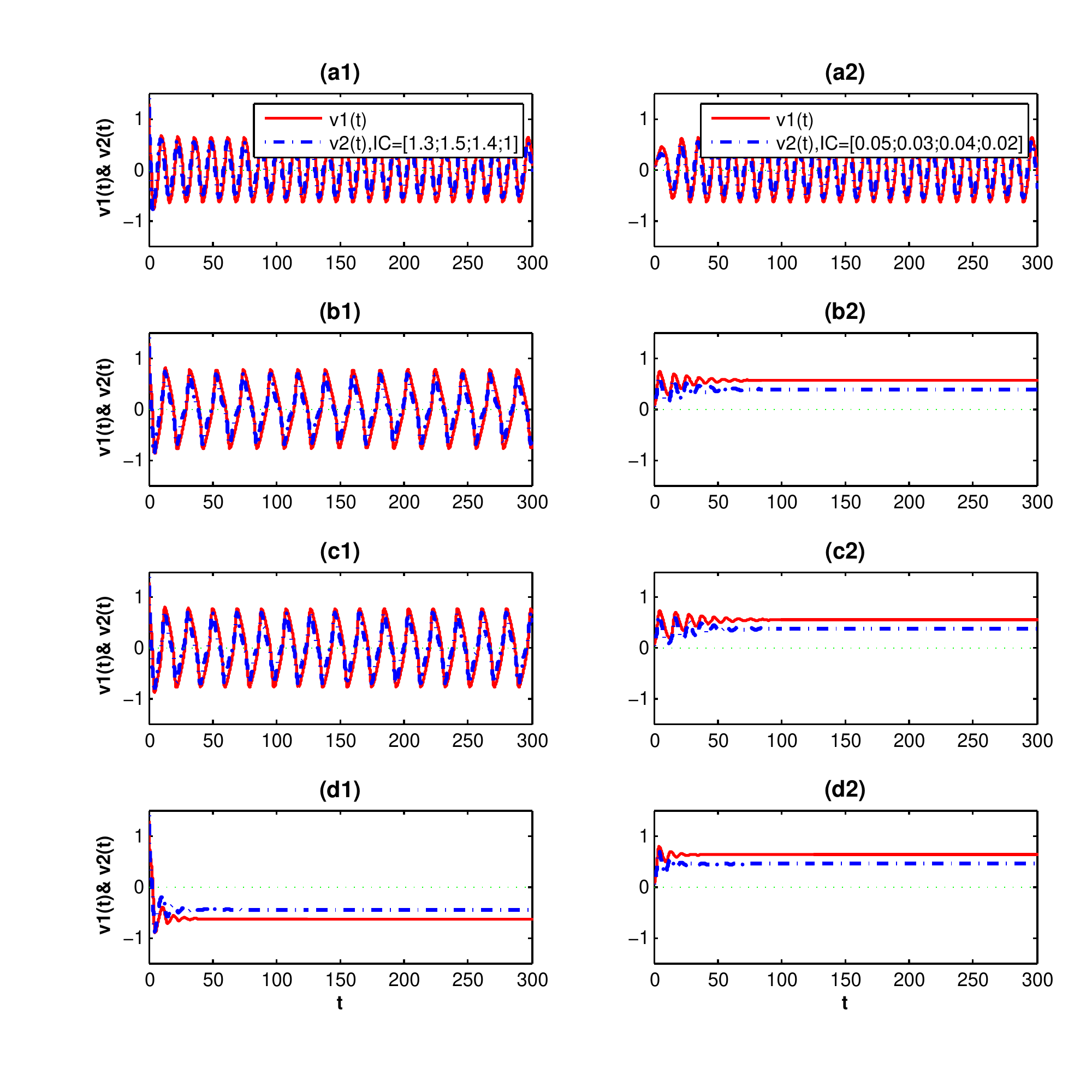}
\caption{\footnotesize Left and right figures show time series of $v_{1}(t)$ (solid red curve) and $v_{2}(t)$ (dashed blue curve) with initial conditions $IC=[ 1.3;1.5;1.4;1]$ and $IC=[0.05;0.03;0.04;0.2]$ respectively. (a1) and (a2) Periodic spiking for parameters in region A of  Fig. \ref{fig6}. (b1) and (b2) Periodic spiking and rest state for parameters in region B of  Fig. \ref{fig6}. (c1) and (c2) Periodic spiking and rest state for parameters in region C of  Fig. \ref{fig6}.  (d1) and (d2) Rest states for parameters in region D of Fig. \ref{fig6}.}
\label{fig7}
\end{figure}

As stated in the previous sections there is a pitchfork bifurcation and two new unstable rest points appear. We want to study the bifurcations of the nontrivial rest points. We will show that the small time delay can induce new period windows with the coupling strength increasing.
In subsection \ref{sec:instantaneous} we stated that there is an eight-figure saddle homoclinic  bifurcation near the sub-critical Hopf bifurcation, for the system with instantaneous coupling. Therefore, in the system with small delayed coupling, in addition to the branch of Hopf bifurcation for small delays, there exists a delay-induced branch of eight-figure saddle homoclinic bifurcation, see Fig. \ref{fig6}. 
Moreover there is a branch of fold of limit cycles that emanates at $c\simeq1.0721$, close to the Hopf and eight-figure saddle homoclinic branches.\\
In order to fully understand the bifurcation diagram in Fig. \ref{fig6}, we change the coupling parameter $c$ and study the dynamic behavior of the system.\\
(1) For $\tau=0.12$ and $c =0.7$ the bifurcation parameters are located below the sub-critical Hopf branch of the non-trivial rest point, region A in Fig.\ref{fig6}. In this region the system is mono-stable and there exists only a stable periodic solution, Fig. \ref{fig7}(a1,a2).\\
(2) For $\tau=0.12$ and $c = 1$, the bifurcation parameters are located above the sub-critical Hopf branch of the non-trivial rest point and below the homoclinic branch, region B in Fig. \ref{fig6}. In this region, the system has a stable limit cycle and a pair of stable nontrivial rest points. With different initial conditions one can reach the stable limit cycle and start to spike or tend to the nontrivial rest points, Fig. \ref{fig7}(b1,b2).\\
(3) For $\tau=0.12$ and $c=1.08$ the bifurcation parameters are located above the homoclinic branch and below the fold of limit cycles branch, region C in  Fig. \ref{fig6}. In this region the behavior of the neurons doesn't change and is like what stated above, Fig. \ref{fig7}(c1,c2).\\
(4) For $\tau=0.12$ and $c=1.1$ the bifurcation parameters are located above the fold branch of limit cycles, region D in Fig. \ref{fig6}, and only non-trivial rest points are stable states in this region, see Fig. \ref{fig7}(d1,d2).
\section{Purely delay driven dynamics}\label{sec:purely}
In this section, in addition to the bifurcations which was stated in previous sections, we will study purely delay driven bifurcations which have no trace in system with instantaneous coupling. 
\subsection{Bifurcations of the trivial rest point}\label{sec:trivial equilibrium}
\begin{figure}
\flushleft
\includegraphics[width=0.5\textwidth]{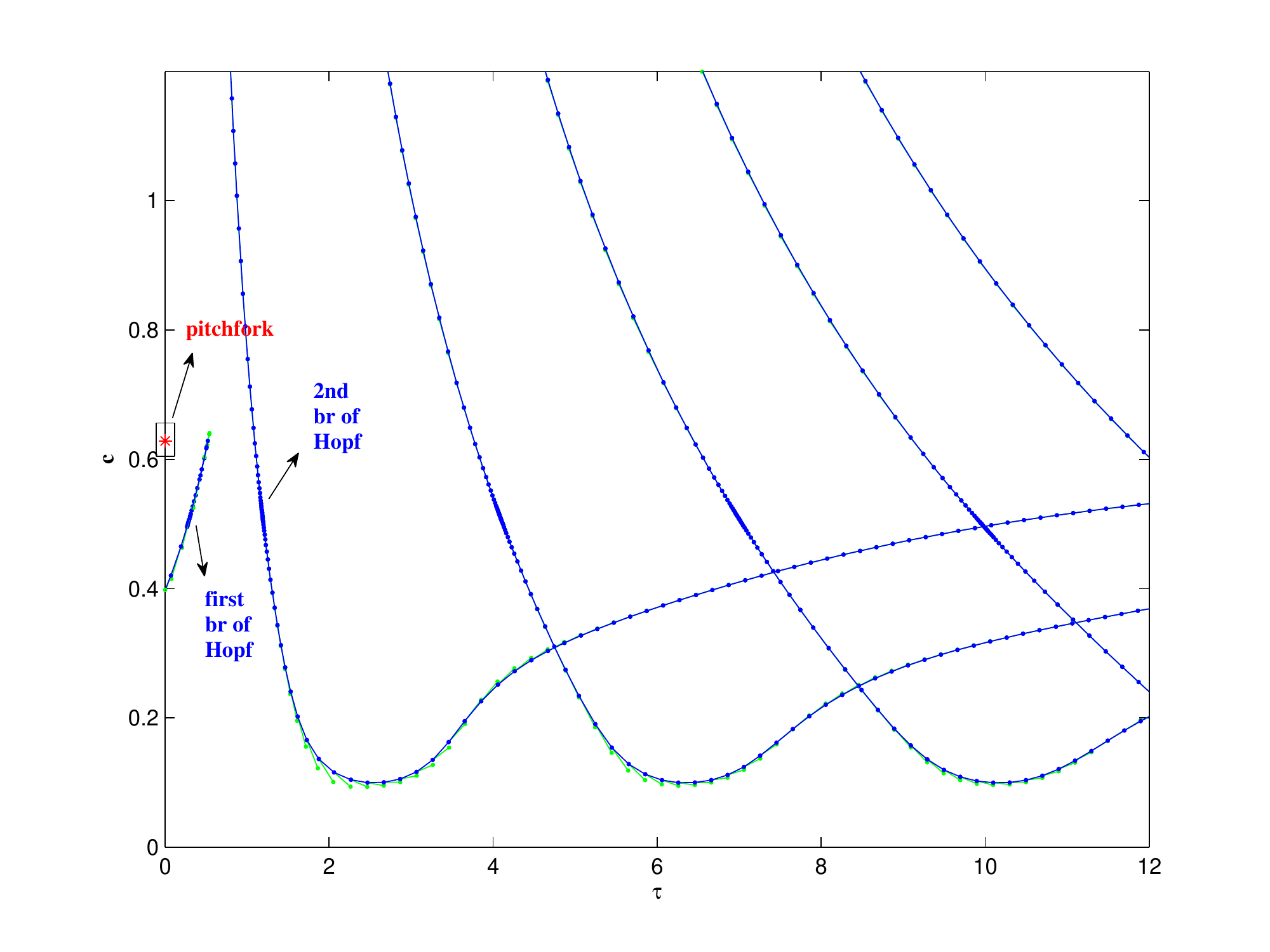}
\caption{\footnotesize Hopf branches of the trivial rest point for $a=0.55, b_{1}=1.128$, and $b_{2}=0.58$. The solid green line is the branch of pitchfork bifurcation. }
\label{fig8}
 \end{figure}

\begin{figure}
\includegraphics[width=0.5\textwidth]{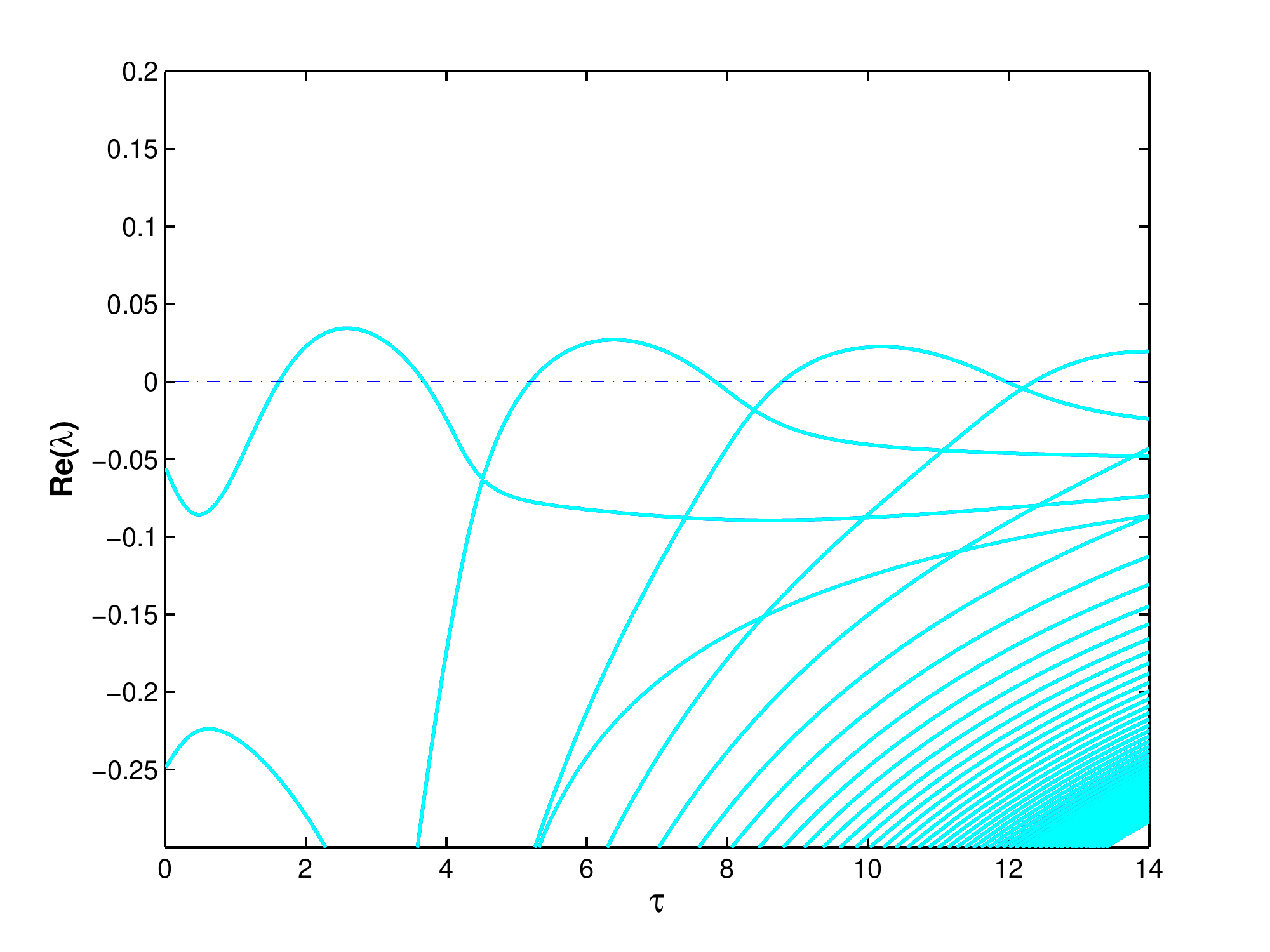} 
\caption{\footnotesize Real parts of eigenvalues for the trivial rest point vs changes of the parameter $\tau$. We set parameters $a=0.55, b_{1}=1.128$, $b_{2}=0.58$, and $c=0.2$.  }
\label{fig9}
\end{figure}

\begin{figure*}
\begin{center}
\includegraphics[width=0.72\textwidth]{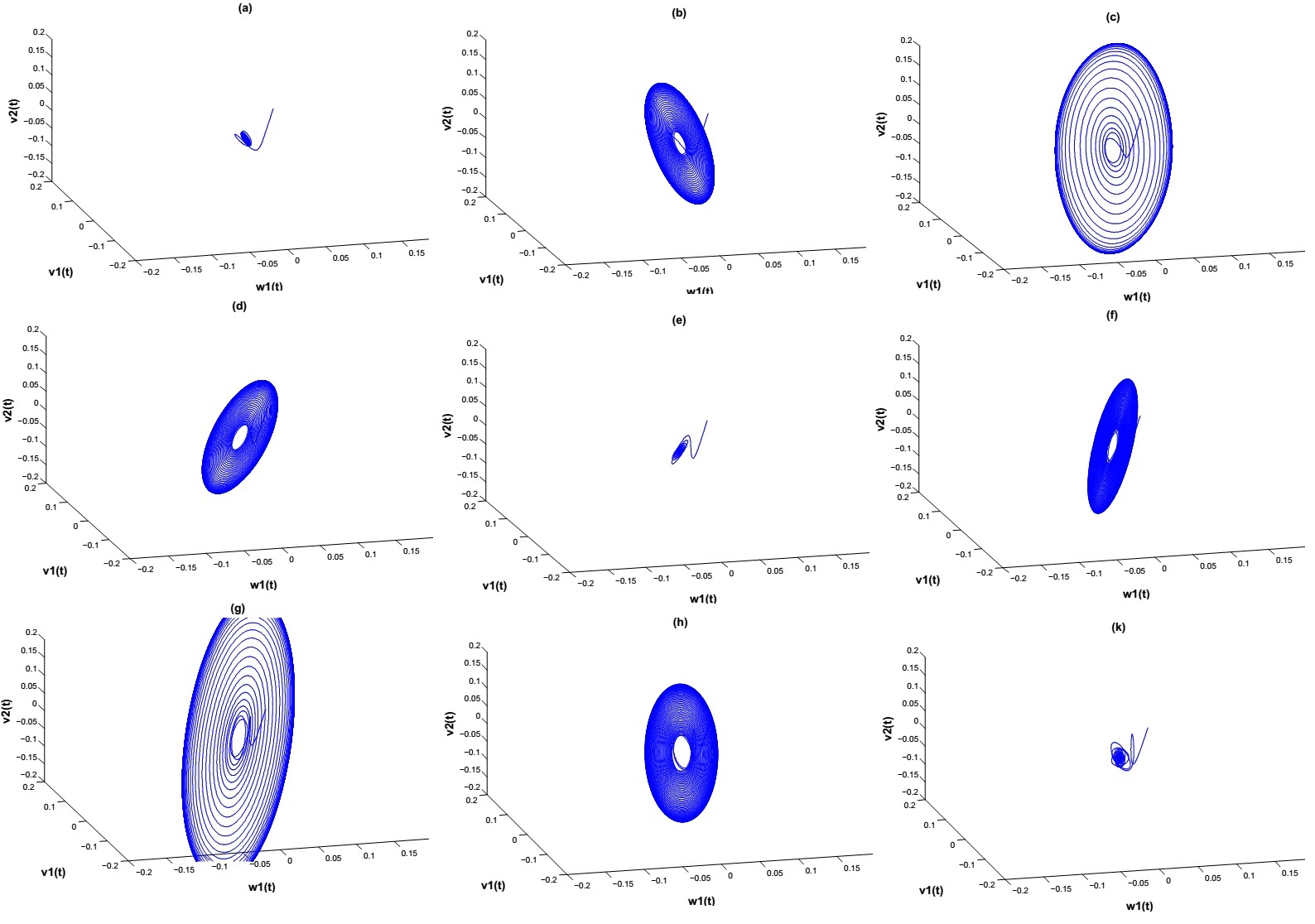}
\end{center}
\caption{\footnotesize Phase portraits of delay-induced oscillations, when $c=0.2$. (a) The origin is a stable rest point for $\tau=1.5<\tau^{(0)}_{2}$. (b) Small periodic oscillation for $\tau=1.8\in(\tau^{(0)}_{2},\tau^{(0)}_{1})$. (c) Large periodic oscillation for $\tau=2.5\in(\tau^{(0)}_{2},\tau^{(0)}_{1})$. (d) The amplitude of the neuron decreases for $\tau=3.5$. (e) The neuron is at rest for $\tau=4\in(\tau^{(0)}_{1},\tau^{(1)}_{2})$. (f) Periodic oscillation for $\tau=5.3>\tau^{(1)}_{2}$. (g) Again, the neuron shows a small periodic oscillation for $\tau=6\in(\tau^{(1)}_{2},\tau^{(1)}_{1})$. (h) The neuron shows a large periodic oscillation. (k) The amplitude of the neuron decreases for $\tau=8.2$.}
\label{fig10}
\end{figure*}

\begin{figure*}
\begin{center}
\includegraphics[width=0.5\textwidth]{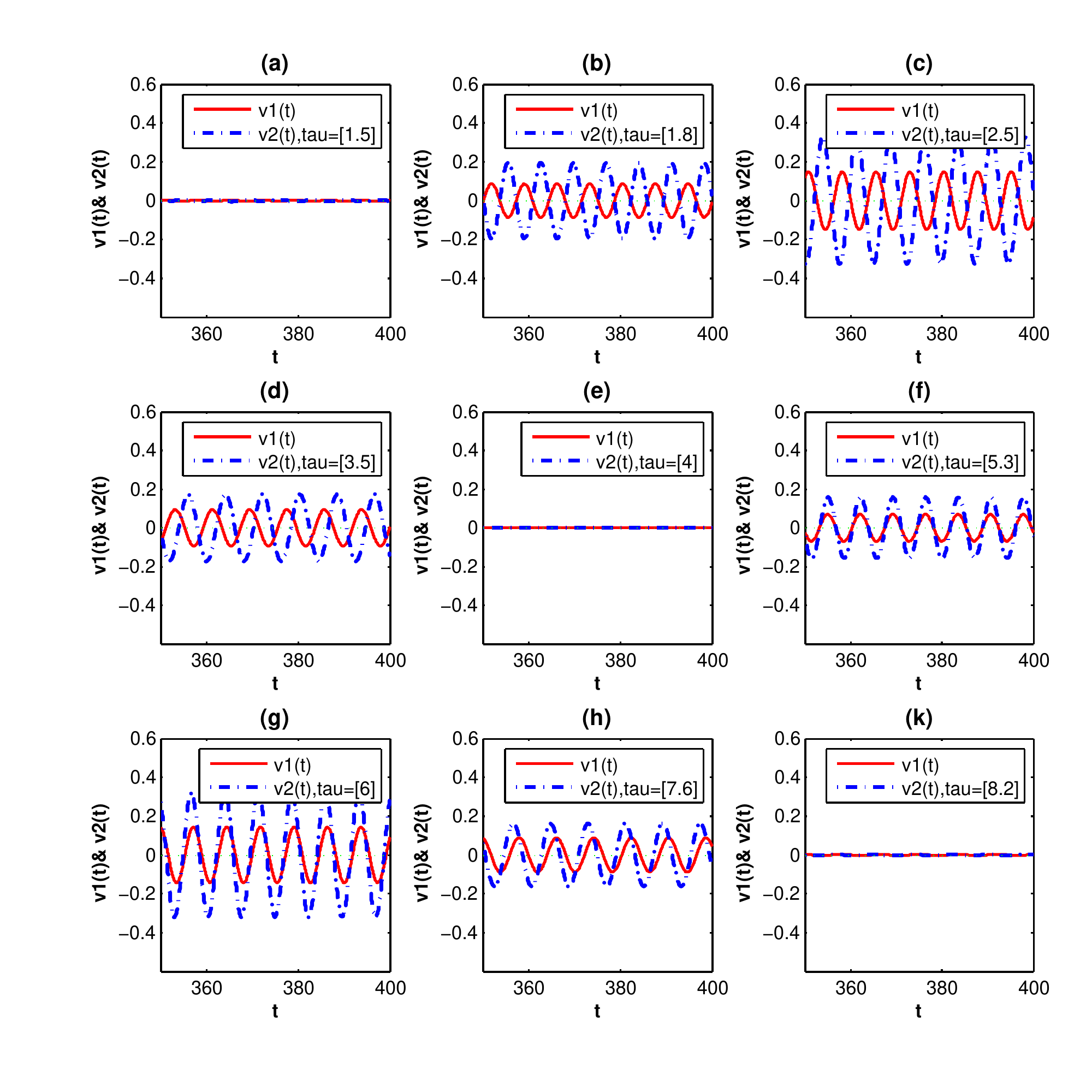}
\end{center}
\caption{\footnotesize Delay-induced oscillations and different modes of oscillation for solutions depicted in Fig. \ref{fig10}. (a) The origin is a stable rest point for $\tau=1.5<\tau^{(0)}_{2}$. (b) The neurons oscillate in an anti-phase manner when $\tau=1.8\in(\tau^{(0)}_{2},\tau^{(0)}_{1})$. (c) The neurons oscillate in an almost anti-phase manner, when $\tau=2.5\in(\tau^{(0)}_{2},\tau^{(0)}_{1})$. (d) Two coupled neurons oscillate  in an almost anti-phase manner, $\tau=3.5$.  (e) The neurons are at rest for $\tau=4\in(\tau^{(0)}_{1},\tau^{(1)}_{2})$. (f) The neuron shows periodic oscillation for $\tau=5.3>\tau^{(1)}_{2}$. (g) Again, the neurons show a small periodic oscillation but synchronous for  $\tau=6\in(\tau^{(1)}_{2},\tau^{(1)}_{1})$. (h) The neurons show large periodic oscillations which are almost synchronized. (k) The amplitude of the neurons decrease for $\tau=8.2$.  }
\label{trivialststrealparts}
\end{figure*}

We start the study of possible bifurcations of the trivial rest point with finding branches of Hopf bifurcation, see Fig. \ref{fig8}.  
The Hopf branches are the roots of equation (\ref{4}), which some of them are shown in Fig. \ref{fig9}.
It is not hard to see that for parameters $c<0.1016$, the conditions of Lemma \ref{stability} hold and the origin is a stable rest point for $\tau\geq0$.\\ 
As the next step, we increase the parameter $c$, to meet the Hopf branches. For instance for $c=0.2$, it is easy to check that the conditions of Lemma \ref{Delta lemma}(i) hold, hence (\ref{4}) has positive roots.\\
We obtain $z^{*}_{1}=0.5739$ and $z^{*}_{2}=0.7718$, thus,$$\omega_{1}=0.7575 \qquad and \qquad \omega_{2}=0.8785.$$
By (\ref{7}), (\ref{astar}), and (\ref{tau}) we have,
\small
\begin{eqnarray}
  \label{taus}  \nonumber \tau^{(j)}_{1}&=&\frac{5.6096+2j\pi}{2(0.7575)};  \tau^{(j)}_{2}=\frac{2.869+2j\pi}{2(0.8785)},\\  j&=&0,1,2,...
\end{eqnarray}
\normalsize
We should notice that for the above parameters, the conditions of Lemma \ref{stability}(ii) are satisfied, hence all roots of (\ref{2}) have negative real parts, when $\tau \in [0,\tau_{0})$, and the neuron is at rest, see Fig. \ref{fig9} and Fig. \ref{fig10}(a) . By Eq. (\ref{taus}) we have $\tau_{0}=\tau^{(0)}_{2}=1.63$.
To understand the bifurcation diagram in Fig. \ref{fig8}, for $c=0.2$ we want to change the parameter $\tau$ and study the dynamic behavior of the system. 
By equations (\ref{taus}) the first Hopf bifurcation occurs for $\tau^{(0)}_{2}$, see Fig. \ref{fig8} and Fig. \ref{fig9}. Actually, it is a super-critical Hopf bifurcation, thus the neuron's stability switches from resting state to periodic spiking, and a small limit cycle branches from the rest point,  Fig. \ref{fig10}(b). Let $T_{2}=\frac{2\pi}{\omega_{2}}$, the resulting limit cycle is  $T$-periodic, with $T$ being near $T_{2}=\frac{2\pi}{\omega_{2}}$. Therefore, the period of the small limit
cycle attractor appearing via super-critical Hopf bifurcation is finite and positive. Thus,
the frequency of oscillations is nonzero, but their amplitudes are small. Such periodic oscillations of small amplitude are not linked to neuronal firing, but must rather be interpreted as spontaneous sub-threshold oscillations. The amplitude of the oscillation grows with the parameter $\tau$ up to a limit,  Fig. \ref{fig10}(c), and after that it decreases, Fig. \ref{fig10}(d), and when the parameter reaches $\tau=\tau^{(0)}_{1}=3.7$, the super-critical Hopf bifurcation occurs and the limit cycle disappears, Fig. \ref{fig10}(e), hence again the trivial rest point becomes stable. We should notice that in the stated region the period of the periodic solution changes and near the parameter $\tau=\tau^{(0)}_{1}=3.7$ reaches $T_{1}$, which $T_{1}=\frac{2\pi}{\omega_{1}}$. By increasing the parameter $\tau$ the trivial rest point remains stable for $\tau \in (3.7,5.2)$. When $\tau=\tau^{(1)}_{2}=5.2$, once more again the super-critical Hopf bifurcation occurs and the behavior of the system changes qualitatively from silence mode to $T_{2}$-periodic spiking, Fig. \ref{fig10}(f). By increasing $\tau$ the above scenario repeats and transitions from silence to firing and vice versa occur frequently until $\tau^{(3)}_{2}=12.36$, Fig. \ref{fig10}(g,h,k). For parameters $\tau>\tau^{(3)}_{2}$, the trivial rest point always has a pair of eigenvalues with positive real part, so the trivial rest point remains unstable and the neurons spike. In this situation, Hopf bifurcations cause the occurrence of a new pair of purely imaginary eigenvalues, Fig. \ref{fig9}.  
Moreover, there is an interesting observation about almost anti-phase and almost synchronized solutions. Actually, almost synchronized and almost anti-phase activities of the coupled neurons can be achieved in some parameter ranges related to their bifurcation transition. If we consider periodic solutions, we can see almost anti-phase activities for $\tau\in(\tau^{(0)}_{2},\tau^{(0)}_{1})$, in which the period changes from $T_{2}$ to $T_{1}$, and almost synchronized activities for $\tau\in(\tau^{(1)}_{2},\tau^{(1)}_{1})$, in which similarly the period changes from $T_{2}$ to $T_{1}$, Fig. \ref{trivialststrealparts}. For an explanation we should consider the eigenspace of the linearized system. The eigenfunction corresponding to the eigenvalue $i\omega_{2}$ is $U^{(j)}(\theta)=e^{i\omega_{2}\theta}(u_{1},u_{2},u_{3},u_{4})$, where $u_{1}=u_{2}(b_{1}+i\omega_{2})$, 
$u_{3}=u_{2}e^{i\omega_{2}\tau^{(j)}_{2}}\frac{1-(a-i\omega_{2})(b_{1}+i\omega_{2})}{c}$, and
$u_{4}=u_{2}e^{i\omega_{2}\tau^{(j)}_{2}}\frac{1-(a-i\omega_{2})(b_{1}+i\omega_{2})}{c(b_{2}+i\omega_{2})}$. It is easy to check that when $j$ is an even number (including zero), $e^{i\omega_{2}\tau^{(j)}_{2}}= \cos(\frac{\omega_{2}}{2}+j\pi) + i\sin(\frac{\omega_{2}}{2}+j\pi) $ lies  in the first quadrant, so if $u_{3}$ and $u_{4}$ lie in the first quadrant, $u_{1}$ and $u_{2}$ will lie in the fourth quadrant and vice versa. Therefore, two neurons have almost anti-phase oscillations. When $j$ is an odd number, $e^{i\omega_{2}\tau^{(j)}_{2}}= \cos(\frac{\omega_{2}}{2}+j\pi) + i\sin(\frac{\omega_{2}}{2}+j\pi) $  lies  in the forth quadrant. Thus, $u_{1}$, $u_{2}$, $u_{3}$, and $u_{4}$ lie in the same quadrant. Therefore, two neurons have almost synchronized oscillations.
As a result, it is seen that time delay can influence dynamical behavior of coupled neurons, such as changing between altering almost synchronized or almost anti-phase oscillations. The stated phenomenon is a consequence of double Hopf bifurcation which is of codimension-two \cite{guo2013bifurcation}.
\subsection{Double Hopf bifurcation}
\begin{figure*}
\begin{center}
\includegraphics[width=0.7\textwidth]{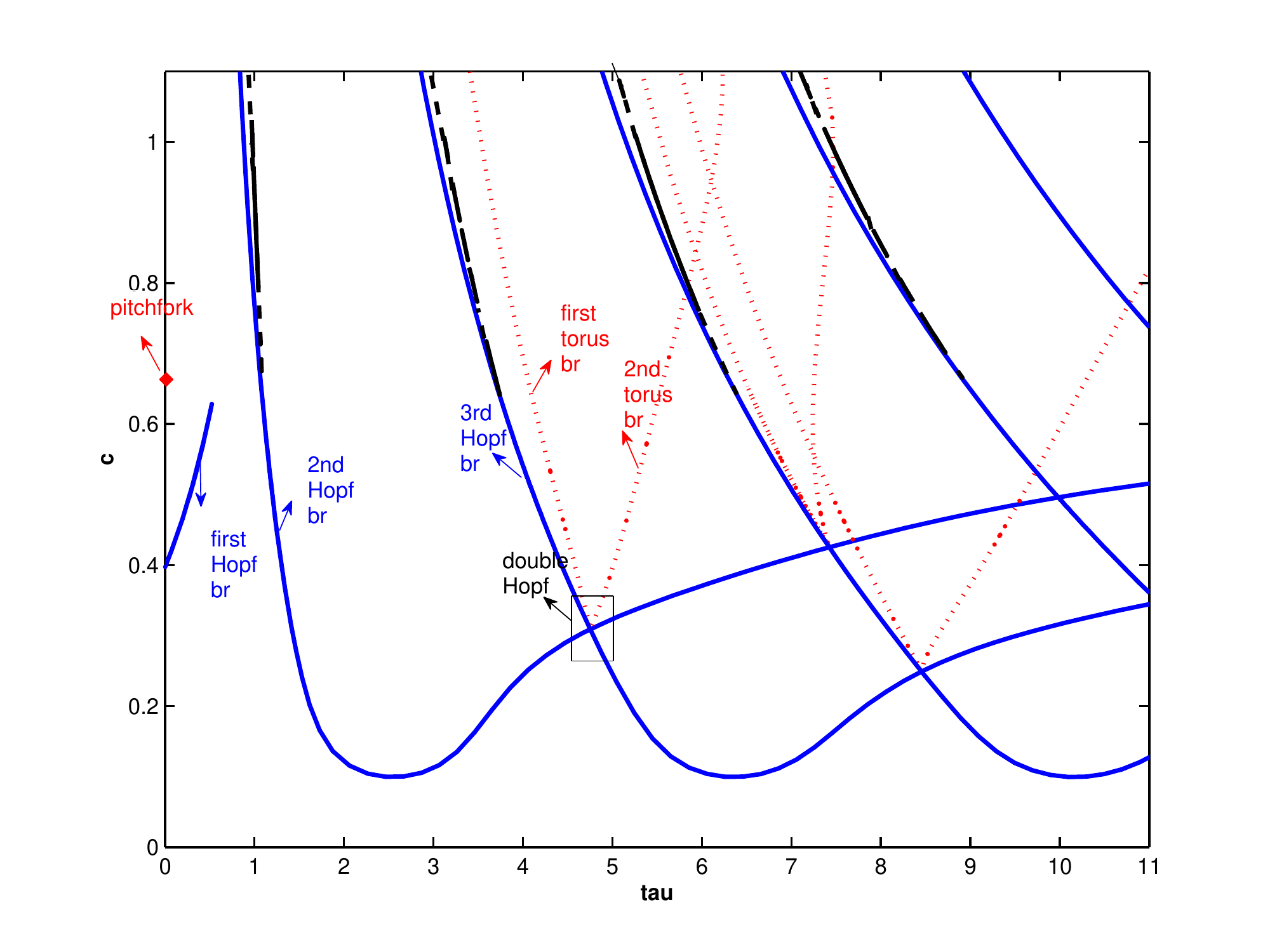}
\caption{\footnotesize Branches of Hopf (blue solid lines), torus (red dotted lines), and pitchfork cycle (black dashed lines) bifurcations of the trivial rest point.  The solid green line is the branch of pitchfork bifurcation. }
\label{torus and Hopf}
\end{center}
\end{figure*}

\begin{figure}
\includegraphics[width=0.5\textwidth]{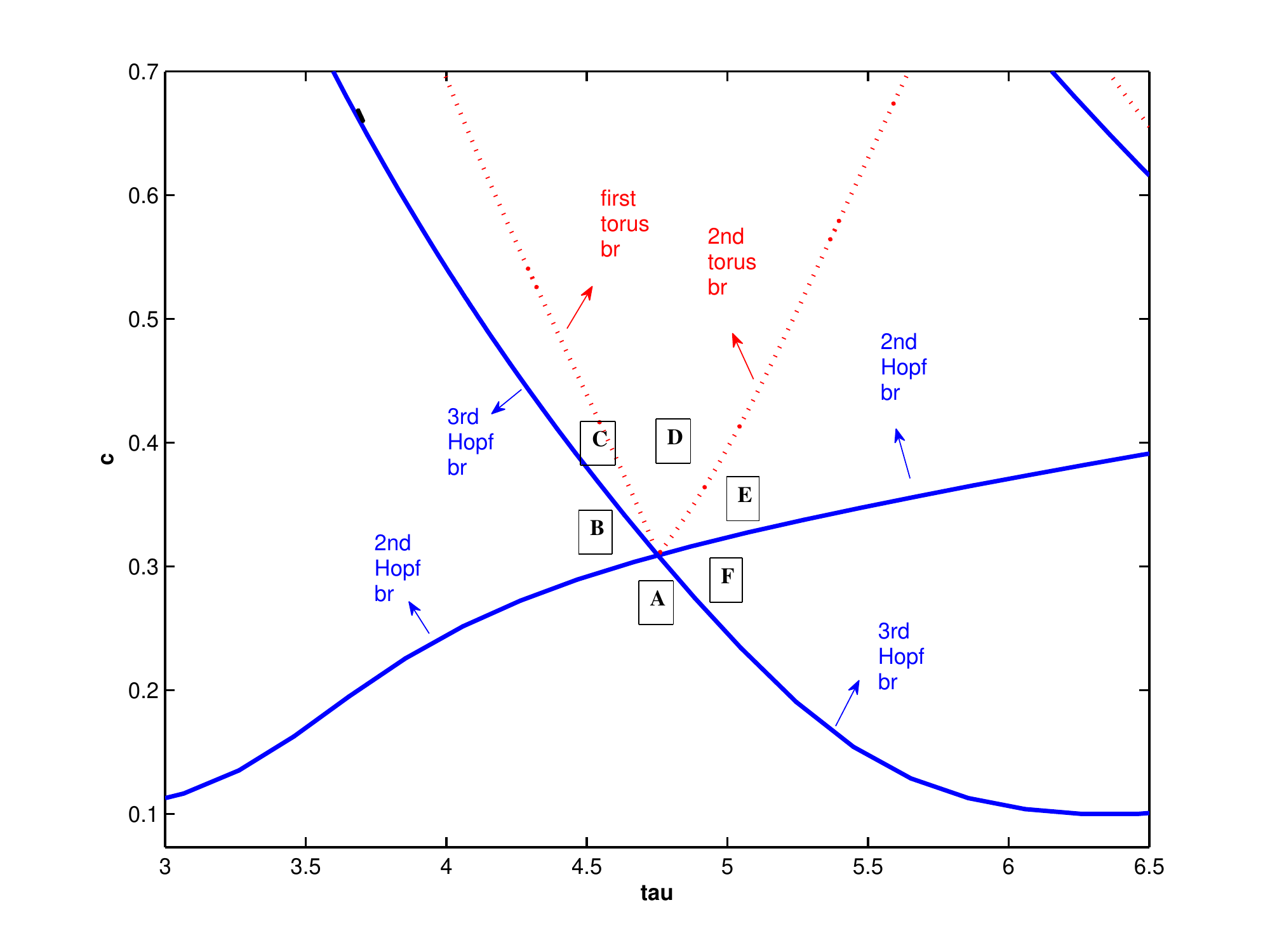}
\caption{ \footnotesize The different regions around a double-Hopf bifurcation of the non-trivial rest point which is related to second and third Hopf branches.}
\label{DH}
\end{figure}

\begin{figure}
\includegraphics[width=0.35\textwidth]{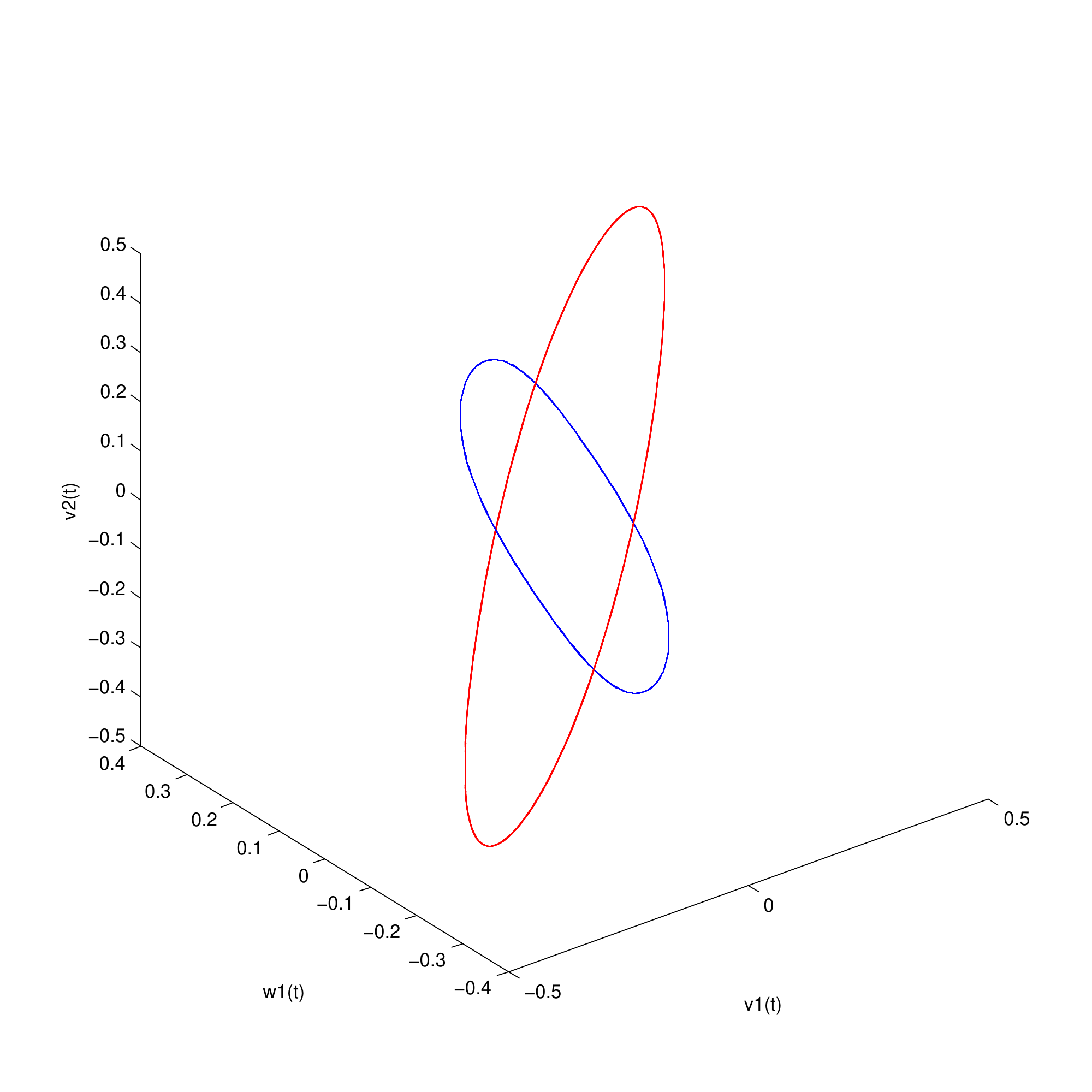} 
\caption{\footnotesize Two different limit cycles with different frequencies which exist when the parameters are in region D. One of them is correspondent to synchronized solutions and another one is correspondent to anti-phase solutions, Fig. \ref{resonance}.}
\label{two stable}
\end{figure}

\begin{figure}
\includegraphics[width=0.4\textwidth]{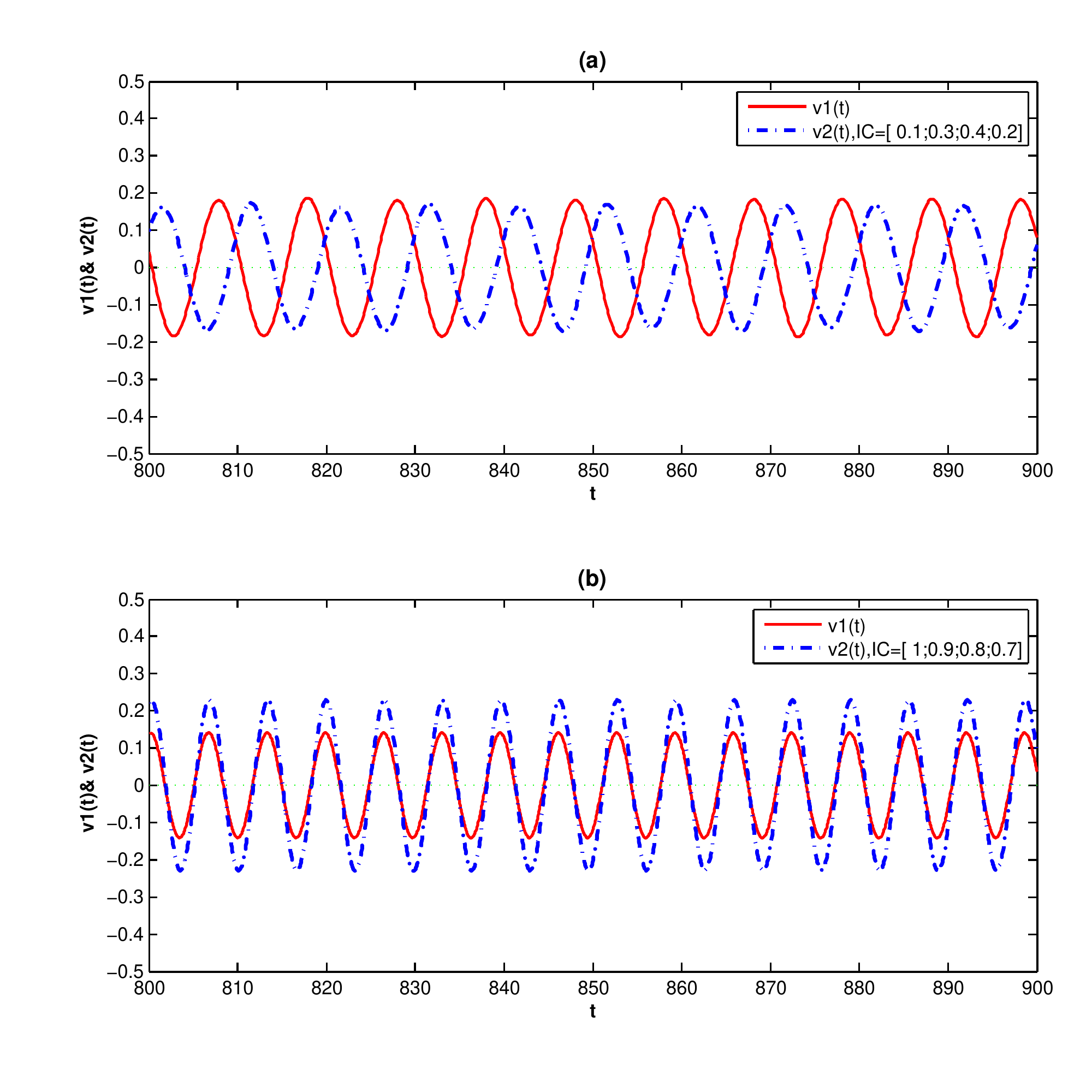} 
\caption{\footnotesize Different modes of spiking due to different initial conditions, here $ c = 0.325$ and $\tau = 4.7756$.  (a) Almost anti-phase oscillations when IC=$ (0.1,0.3,0.4,0.2)$. (b) Almost synchronized oscillations when IC=$(1,0.9,0.8,0.7)$.  }
\label{resonance}
\end{figure}

\begin{figure}
\includegraphics[width=0.4\textwidth]{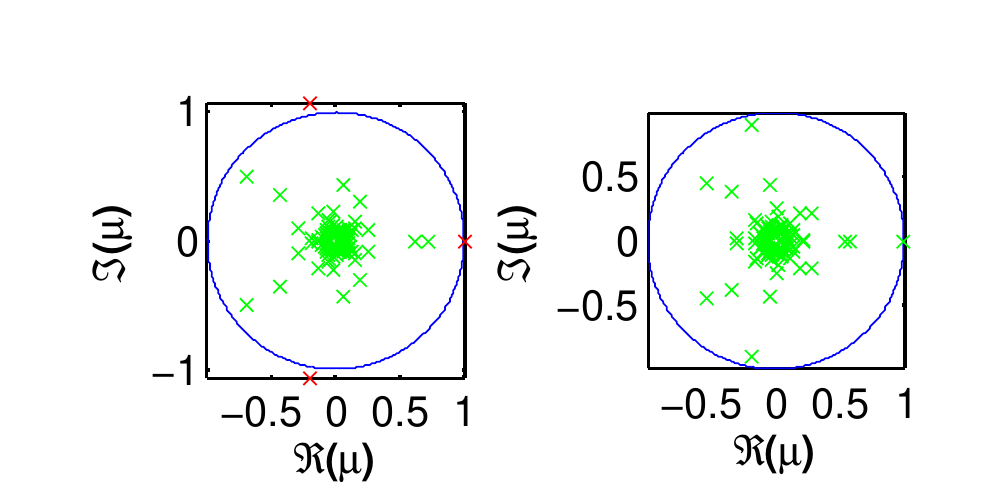} 
\caption{\footnotesize (Left) Floquet multipliers before first torus bifurcation. (Right) Floquet multipliers after first torus bifurcation. }
\label{from third hopf}
\end{figure}

In the previous sections we reviewed Hopf bifurcations of the trivial rest point and the multiple stability regions which are formed in the $(\tau, c)$ plane, and imply multi-time switches between the stable and unstable states of the zero solution. In this section we want to study other possible bifurcations of the trivial rest point, Fig. \ref{torus and Hopf}.
 There are some points in the bifurcation diagram of the trivial rest point, Fig. \ref{torus and Hopf}, which can be obtained from the intersection of the two branches of Hopf bifurcation.
It is possible that the characteristic equation at a rest point has eigenvalues with strict nonzero real parts except two pairs of purely imaginary eigenvalues; in $(\tau, c)$ plane it happens when two branches of Hopf bifurcation denoted by $\tau_{1}$ and $\tau_{2}$ correspond to $\omega_{1}$ and $\omega_{2}$ respectively, cross each other in a point $(\tau_{dh}, c_{dh})$. This situation is called a double-Hopf bifurcation \cite{guo2013bifurcation}. For constants $(\tau_{dh}, c_{dh})$, let $\pm i \omega_{1}$ and $\pm i \omega_{2}$ be the two pairs of purely imaginary eigenvalues of the characteristic equation and all other eigenvalues have non zero real parts. The existence of double-Hopf bifurcation leads to different modes of spiking in different regions around the double-Hopf point, regions C,D, and E in Fig. \ref{DH}.\\
It follows from the previous paragraph that different values of delay and coupling parameters of the system (\ref{1}), can induce two different values of frequencies corresponding to the second and third Hopf bifurcation. The existence of a double-Hopf bifurcation can be detected by the local stability analysis. In this section, we represent the typical dynamical behaviors and the bifurcation sets in the neighborhood of one of the double-Hopf bifurcation points of Fig. \ref{torus and Hopf}.
First, time delay and coupling strength are chosen from the region A of Fig. \ref{DH}. In this region the trivial rest point is locally asymptotically stable. By increasing the parameter $c$, the trivial rest point loses its stability in result of the super-critical Hopf bifurcation, where the parameters encounter the second Hopf bifurcation branch, and enter region B of Fig. \ref{DH}. The system periodically oscillates with the $\omega_{1}$ frequency. On the other hand, when the parameter $\tau$ is increased from region A to F of Fig. \ref{DH}, the parameters will encounter the third Hopf bifurcation branch related to $\omega_{2}$ frequency. This implies that the trivial rest point loses its stability through a super-critical Hopf bifurcation and the system exhibits the oscillatory solution with  $\omega_{2}$ frequency. If we start from region B with a fixed parameter $c$ and increase the parameter $\tau$, the parameters encounter the third super-critical Hopf bifurcation branch and enter region C. Since the trivial rest point is unstable, an unstable limit cycle appears. In this situation, the system is mono-stable and exhibits the oscillatory solution with  $\omega_{1}$ frequency. By further increasing the the parameter $\tau$ from region C to D, the unstable periodic solution with $\omega_{2}$ frequency becomes stable through the sub-critical torus bifurcation. It suggests that the system behavior changes to co-existence of two different periodic solutions with different frequencies. Actually for parameters in region D there are two different stable limit cycles as depicted in Fig. \ref{two stable}, and the system is bi-stable. 
We should notice that in this region the behavior of the system can be different when we start with different initial conditions and two neurons without changing the parameters can show oscillatory solutions with two different frequencies, and different types of synchronization, Fig. \ref{resonance}. 
We should emphasis that there is an unstable torus in region C. Torus bifurcation occurs when a complex-conjugate pair of Floquet multipliers crosses the unit circle \cite{kuznetsov2013elements}. In the above stated torus, a complex-conjugate pair of Floquet multipliers moves from the outside to the inside of the unit circle, Fig.\ref{from third hopf}.
 \begin{figure}
\includegraphics[width=0.4\textwidth]{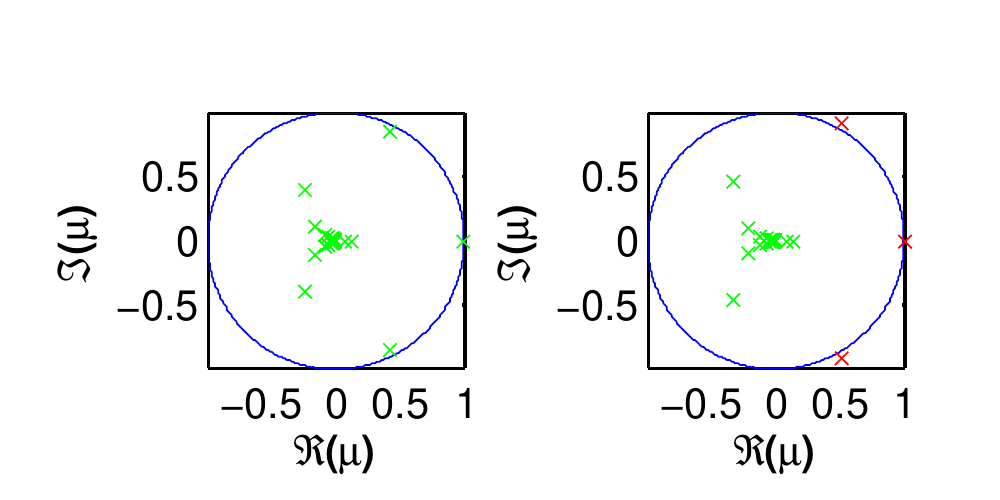}
\caption{\footnotesize (Left) Floquet multipliers before second torus bifurcation. (Right) Floquet multipliers after second torus bifurcation. }
\label{from second hopf}
\end{figure}
When the time delay is increased from region
D to E of Fig. \ref{DH}, the stable periodic solution which appeared through the second branch of Hopf bifurcation loses its stability through the sub-critical torus bifurcation. Hence, in region D the system becomes mono-stable. For the stated torus a complex-conjugate pair of Floquet multipliers moves from inside to the outside of the unit circle, Fig. \ref{from second hopf}.\\
As we showed in the figures, Hopf bifurcation gives rise to the periodical activity, while double-Hopf bifurcation results in two types of oscillations with different frequencies and different synchrony, and a small change in the time delay can result in different patterns of activity. This is particularly useful when using the time delay to control bifurcations and oscillatory solutions.\\
As it is depicted in Fig. \ref{torus and Hopf}, there are some branches of  pitchfork cycle bifurcation. We will explain their roles on the system's dynamic in the next section.

\subsection{Bifurcations of the non-trivial rest points}
\begin{figure}
\includegraphics[width=0.55\textwidth]{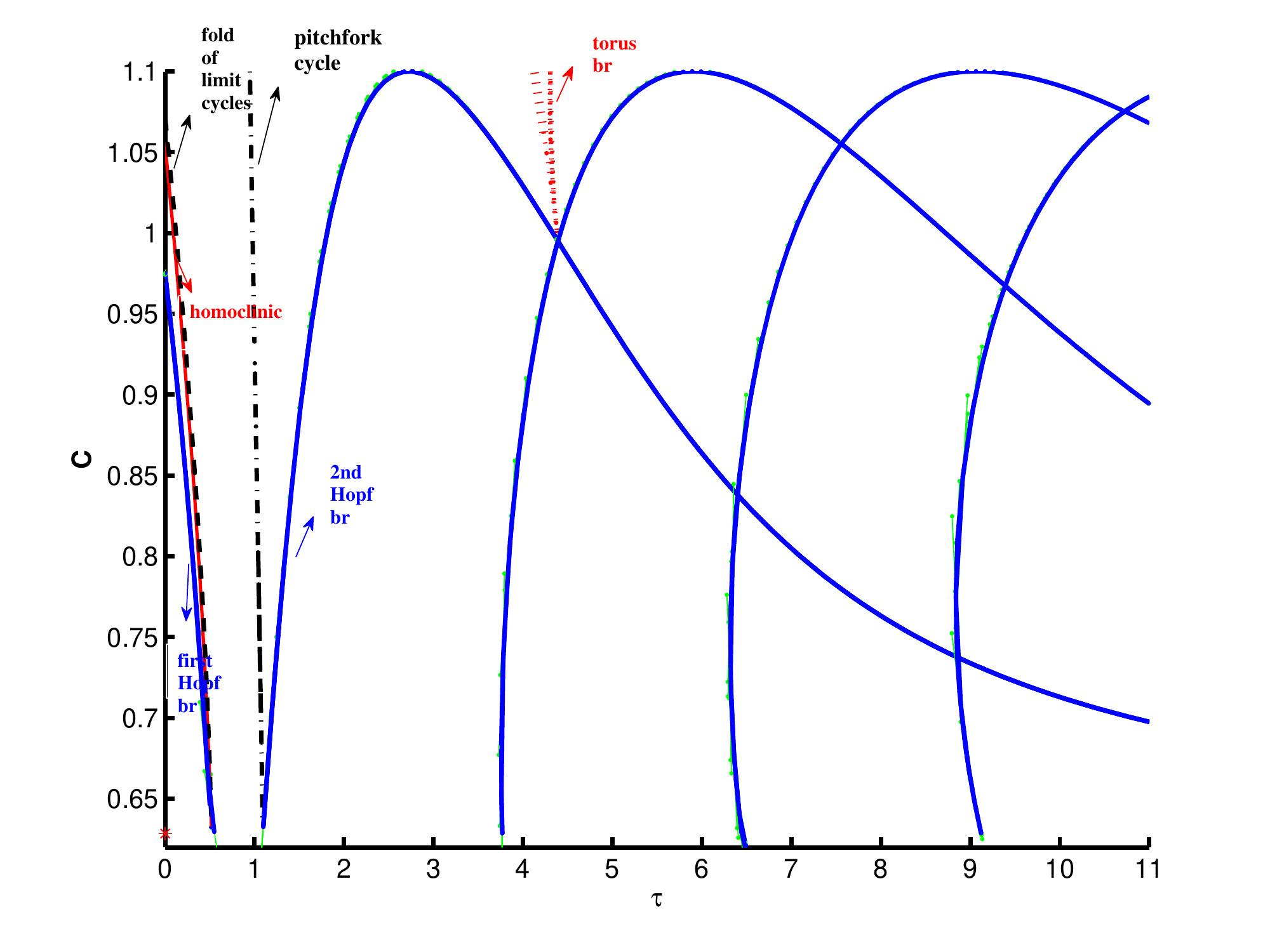}
\caption{\footnotesize Bifurcations of the non-trivial rest points. The solid blue lines are Hopf bifurcation branches. The dash-dot black line is the pitchfork cycle bifurcation. The dashed black line is fold bifurcation of limit cycles. The dotted red lines are torus bifurcation branches. The solid red line is fold of limit cycles branch. }
\label{fig18}
\end{figure}

\begin{figure*}
\begin{center}
\includegraphics[width=0.7\textwidth]{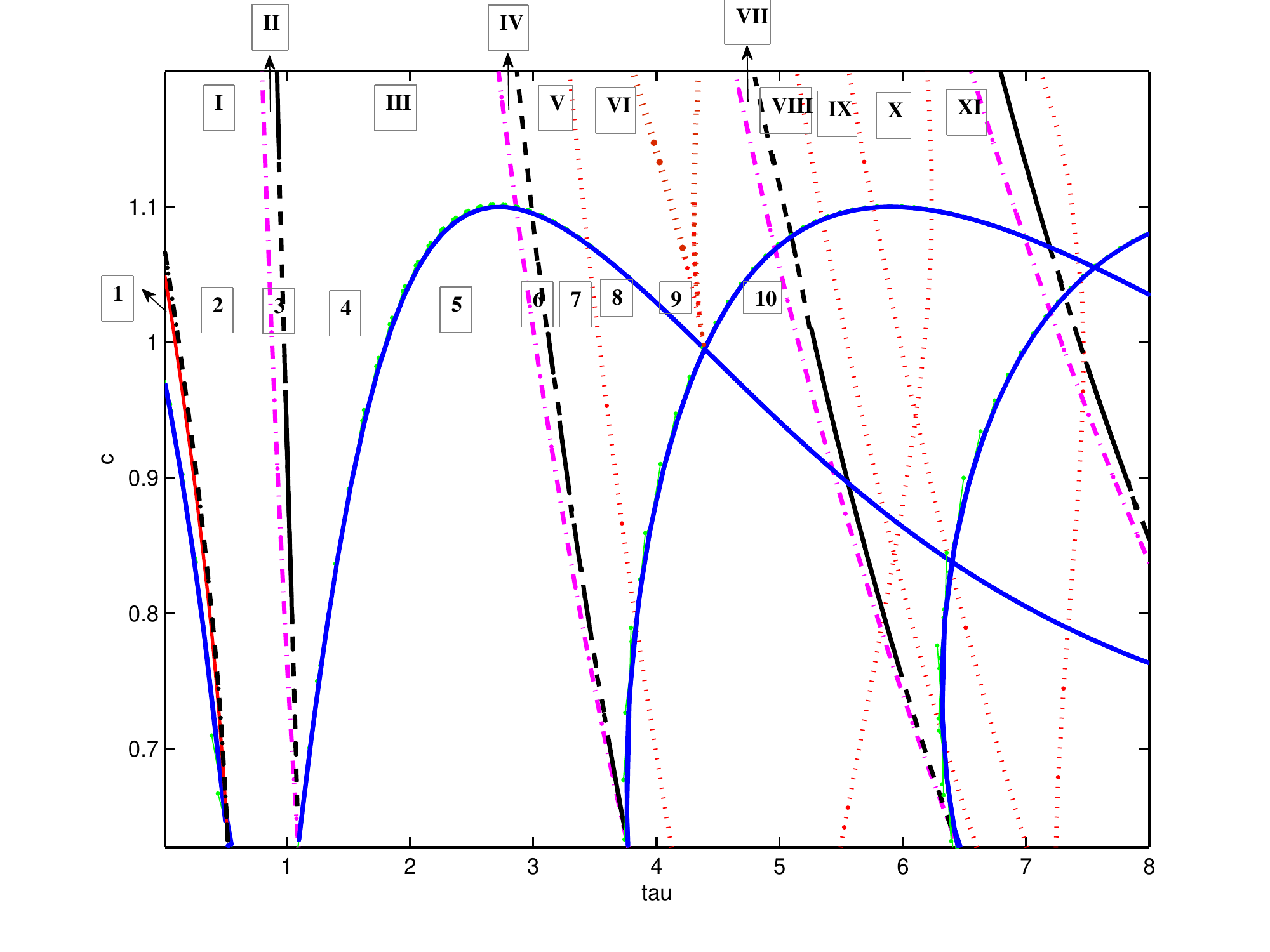}
\end{center}
\caption{\footnotesize Total bifurcation diagram of the system (\ref{1}). Blue solid lines are Hopf branches of non-trivial rest points. Magenta dash-dot lines are Hopf branches of trivial rest point. Black dashed lines are pitchfork bifurcations of limit cycles. Red dotted lines are torus bifurcation branches. The other branches are introduced in Fig. \ref{fig6}.}
\label{new}
\end{figure*}

\begin{figure}
\includegraphics[width=0.5\textwidth]{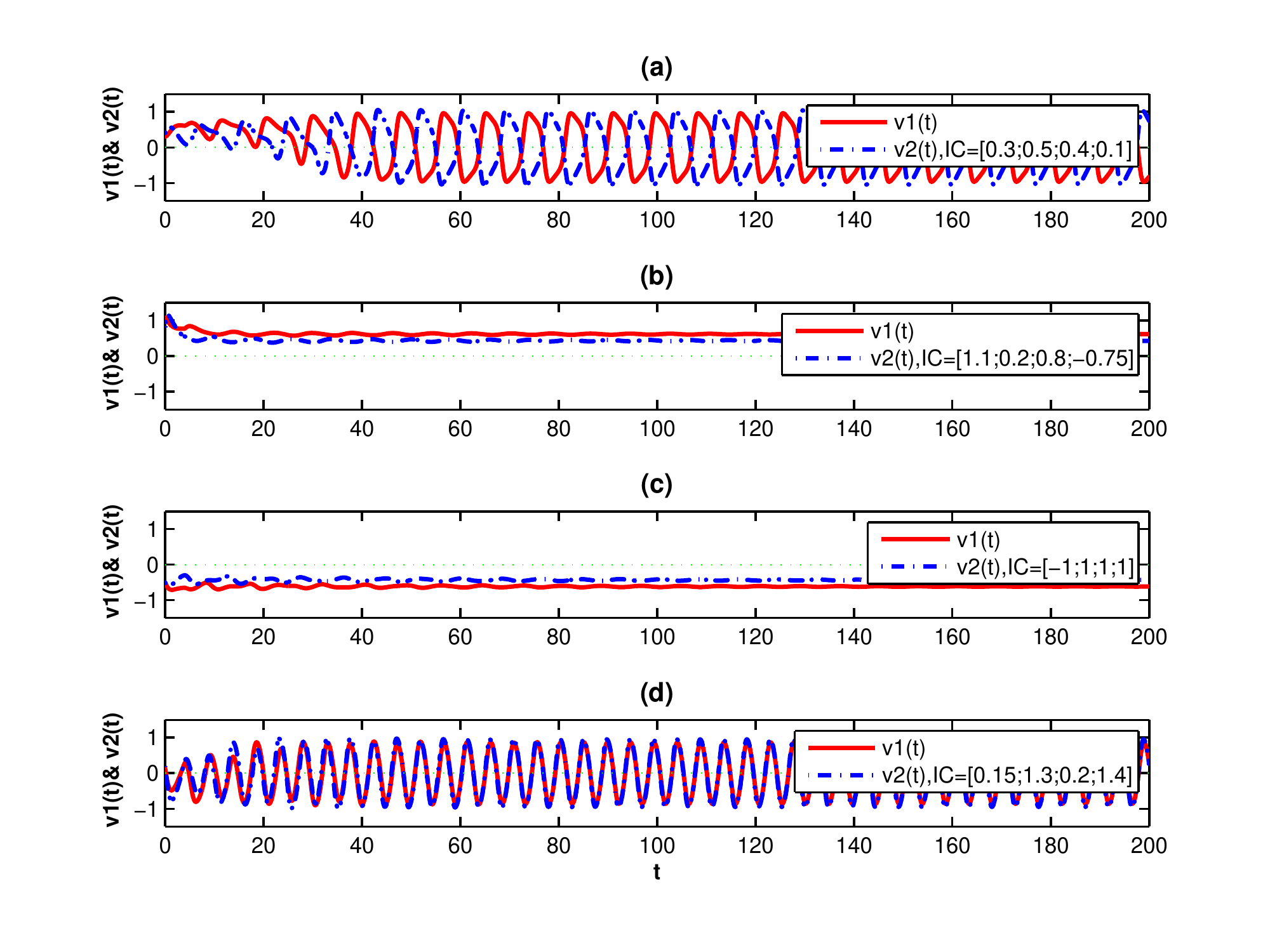}
\caption{\footnotesize Different dynamical activities in region 9 of  Fig. \ref{new}. Whit different initial conditions the system reaches different stable states. Here $c=1.08$ and $\tau=3.9$.  }
\label{new1}
\end{figure}

In section \ref{sec:trivial equilibrium}, according to the bifurcations of the trivial rest point, we showed different stability regions in the $(\tau, c)$ plane, which implies multi-time switches between the stable and unstable states of trivial solution. Similarly we can see such multiple stability regions according to bifurcations of  non-trivial rest points. In this way, we will study the possible bifurcations of the non-trivial rest points too, see Fig. \ref{fig18}.
 In the following, we restrict ourselves to a fixed value of $c$ and study the delay induced changes of the dynamics according to all possible bifurcations of the system. Since in the study of the dynamics of the system (\ref{1}), the bifurcations of the origin and non-trivial rest points are both indispensable, we draw entire bifurcation diagram which were depicted in Fig. \ref{fig18} and Fig. \ref{torus and Hopf} in a single figure, see Fig. \ref{new}.
 For fixed parameter $c=1.03$, we change the parameter $\tau$, and investigate the behavior of the system. When $\tau=0.05$, the parameters are located in region 1 of  Fig. \ref{new} and  the dynamics is the same as in region B in Fig. \ref{fig6}. Therefore, the system has a stable limit cycle related to synchronous oscillations, and a pair of stable non-trivial rest points. By increasing the parameter $\tau$, for example $\tau=0.5$, the parameters are located in region 2 of  Fig. \ref{new}, and  the dynamics is as region D in Fig. \ref{fig6}. Since the parameters are located after the branch of fold of limit cycles, only non-trivial rest points are stable. By a little change of the parameter $\tau$ the second branch of super-critical Hopf bifurcation of the trivial rest point occurs, and the trivial rest point which was unstable gives birth to a small-amplitude unstable limit cycle. Therefore, in region 3 of  Fig. \ref{new} the two non-trivial rest points remain stable.  By increasing the parameter $\tau$ a sub-critical pitchfork bifurcation of limit cycles occurs, the stated unstable limit cycle changes stability, becomes stable and produces two unstable limit cycles. In result, in this region, according to the initial conditions the neurons spike in an anti-phase manner, or remain in the rest states (non-trivial rest points), region 4 in Fig. \ref{new}.  By more changes of $\tau$ the delay-induced limit cycle stay stable and never disappear, but each of the unstable limit cycles goes towards a non-trivial rest point and their amplitude decrease. When $\tau$ reaches the second branch of sub-critical Hopf bifurcation of non-trivial rest points in Fig. \ref{fig18}, the small unstable limit cycles shrinks to stable non-trivial rest points and make them lose stability. Hence, the system becomes mono-stable and the neurons start to spike, region 5 in Fig. \ref{new}.  By further increasing of the parameter $\tau$, the third branch of super-critical Hopf bifurcation for trivial rest point appears,  Fig. \ref{torus and Hopf}, and the trivial rest point, which is unstable of Morse index 3, gives birth to a small-amplitude unstable limit cycle with index 3. Therefore, in region 6 of  Fig. \ref{new}, the system remains mono-stable.  By increasing the parameter $\tau$ again a sub-critical pitchfork bifurcation of limit cycles occurs, and the index-3 unstable limit cycle changes to index-2 unstable limit cycle. Hence, in region 7 of  Fig. \ref{new} the system remains mono-stable. If we increase the parameter $\tau$, a sub-critical torus bifurcation occurs, and the index-2 unstable limit cycle changes stability and becomes stable. The dynamics in region 8 of  Fig. \ref{new}, is the same as dynamics in region D of Fig. \ref{DH}.  Therefore, the system behavior changes to co-existence of two different periodic solutions with different frequencies and different synchrony. Actually, for parameters in region 8 there are two different stable limit cycles, one of them is related to almost synchronous oscillations and the other one is related to almost anti-phase activities. The dynamics in region 9 of  Fig.\ref{new} is really interesting and complicated. By changing parameters from region 8 to 9, a sub-critical Hopf bifurcation for non-trivial rest points occurs and makes the unstable non-trivial rest points stable. Hence, in region 9, there are  two different periodic solutions with different frequencies, and also a pair of stable rest points, see Fig. \ref{new1}. Related to double-Hopf bifurcation point which is stated below the region 9, there are branches of super-critical torus bifurcations. Since the Hopf bifurcations of the non-trivial rest points are sub-critical, the limit cycles which are related to these bifurcations and also the created torus are unstable. If we increase $\tau$ from region 9 to 10 again a  sub-critical Hopf bifurcation occurs for non-trivial rest points and makes them unstable.  Hence in region 10 of  Fig. \ref{new} there are two different stable limit cycles, and rest points are unstable.
\subsection{Strong coupling}
For the system (\ref{1}) with instantaneous coupling ($\tau=0$) the  characteristic equation  at the trivial rest point is
 \begin{eqnarray} \label{10}
  P(c,\tau)&=& \lambda^{4}+ A \lambda^{3}+ B \lambda^{2}+ C \lambda + D =0, 
\end{eqnarray}
where $ A=b_{1}+b_{2}-2a $, $ B=b_{1}b_{2}-2a(b_{1}+b_{2})+a^{2}-c^{2}+2 $, $ C= (a^{2}+1)(b_{1}+b_{2})-2ab_{1}b_{2}-2a- c^{2}(b_{1}+b_{2})$, and $ D= (a^{2}-c^{2})b_{1}b_{2}-a(b_{1}+b_{2})+1$. According to Lemma \ref{lemma1}, if $D<0$, then the characteristic equation (\ref{10}) has at least one positive root. It is easy to check that for chosen values of the parameter in our paper, when $c>0.6265$, $D$ is negative and the characteristic equation has at least a positive root. Indeed, for $c>0.6265$ the trivial rest point is unstable of index 1. Moreover, in the system with delayed coupling (\ref{1}), when pitchfork bifurcation occurs the trivial rest point becomes unstable of index k, where k is an odd number. Since the delay driven bifurcations such as Hopf bifurcation change the rest point of index k to a rest point of index $k\pm2$, the index of the rest point remains odd, and the trivial rest point remains unstable. Although the trivial rest point is unstable for all the parameter space in Fig. \ref{new}, the influence of it's Hopf bifurcations on the dynamics isn't negligible. 
\begin{figure}
\includegraphics[width=0.4\textwidth]{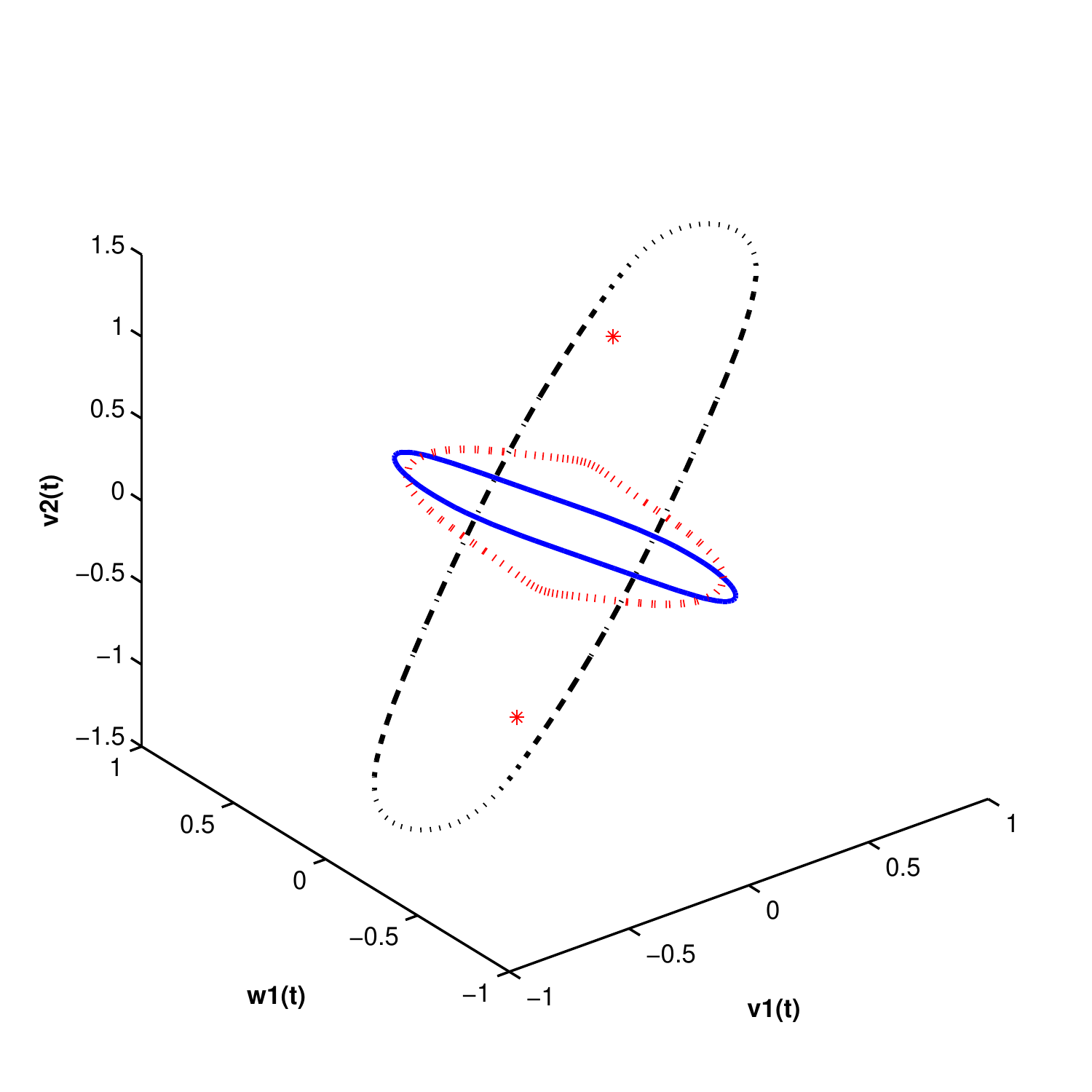}
\caption{\footnotesize  Different limit cycles with different frequencies which exist in region X of  Fig. \ref{new}. The dash-dot black line is correspondent to synchronized solutions and the others are correspondent to anti-phase solutions.  }
\label{fig20}
\end{figure}
\begin{figure}
\includegraphics[width=0.35\textwidth]{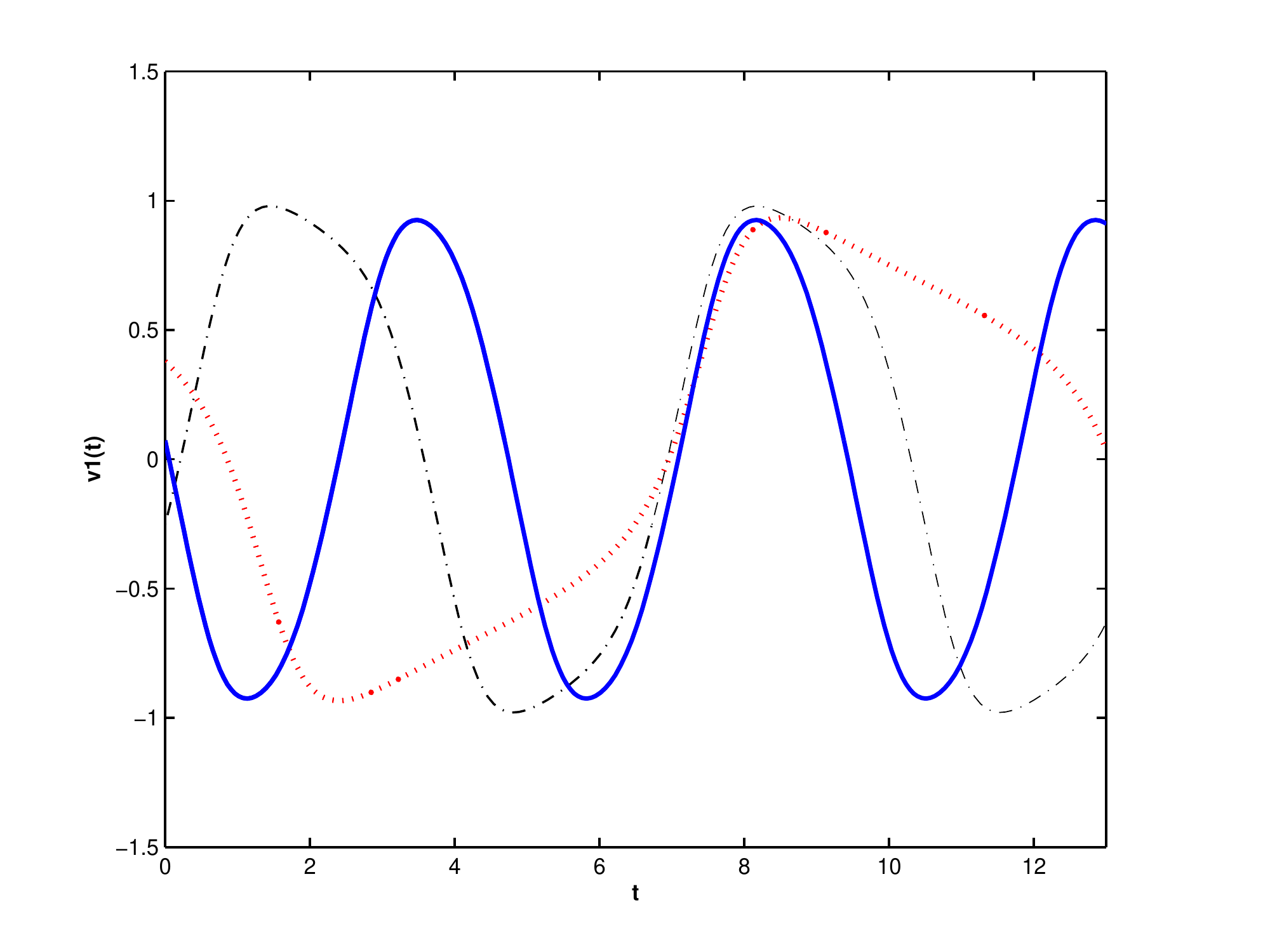}
\caption{\footnotesize  The dash-dot black line is correspondent to synchronized solutions, the solid blue line is correspondent to  anti-phase solutions with highter frequency, and the dotted red line is correspondent to anti-phase solutions with lower frequency.  }
\label{fig21}
\end{figure}
Our interest is to determine if there are some stable solutions, and synchrony properties of oscillatory activities of the neurons.
In this way, for fixed parameter $c=1.2$, we change the parameter $\tau$, and investigate the behavior of the system. When the parameters are located in region I of  Fig. \ref{new}, the dynamics is the same as region 2 of Fig. \ref{new}. Therefore, only  the two non-trivial rest points are stable and they remain stable for all values of parameter $\tau$ in parameter space of  Fig. \ref{new}. By further change of the parameter $\tau$ the second  super-critical Hopf bifurcation of the trivial rest point occurs, and the trivial rest point which was unstable of index 1 gives birth to a small-amplitude unstable limit cycle of index 1. Therefore in region II of  Fig. \ref{new}, the two non-trivial rest points are stable, and there exists an unstable limit cycle.  By a little change of the parameter $\tau$ a sub-critical pitchfork bifurcation of limit cycles occurs, the stated unstable limit cycle changes stability, becomes stable and produces two unstable limit cycles. In result, in  region III of Fig. \ref{new}, according to the initial condition the neurons spike in an almost anti-phase manner, or remain in the rest state (non-trivial rest points).  By further increasing of the parameter $\tau$, the third branch of super-critical Hopf bifurcation of the trivial rest point appears,  Fig. \ref{torus and Hopf}, and the trivial rest point, which is unstable of index 3, gives birth to a small-amplitude unstable limit cycle with index 3. Therefore, the activity of the neurons in region IV of  Fig. \ref{new}, remains similar to those in region III.  By increasing the parameter $\tau$, again a sub-critical pitchfork bifurcation of limit cycles occurs, and the index-3 unstable limit cycle changes to index-2 unstable limit cycle. Hence, in region V of  Fig. \ref{new} the activity of the neurons is the same as region IV. If we increase the parameter $\tau$, a sub-critical torus bifurcation occurs, and the index-2 unstable limit cycle changes stability and becomes stable.  This stable solution is related to almost synchronized oscillations. The dynamics in region VI of  Fig. \ref{new}, is the same as dynamics in region 9, hence there are  two different periodic solutions with different frequencies and different synchrony, and also a pair of stable rest points, see Fig. \ref{new1}. As $\tau$ is increased from region VI to VII,  another branch of super-critical Hopf bifurcation of the trivial rest point appears, and the trivial rest point which is unstable of index 5 becomes unstable of index 7.  This Hopf bifurcation results in the creation
of a  small-amplitude unstable limit cycle with index 5.  By a little change of the parameter $\tau$, a sub-critical pitchfork bifurcation of limit cycles occurs, and the index-5 unstable limit cycle changes to index-4 unstable limit cycle, region VIII in Fig. \ref{new}.  By further change of the parameter $\tau$, a torus bifurcation occurs and the index-4 unstable limit cycle changes to index-2 unstable limit cycle, region IX in Fig. \ref{new}. The activity of the neurons is not affected qualitatively by this bifurcation. When the parameter $\tau$ is increased from region IX to X, the index-2 unstable limit cycle changes stability and becomes stable through a torus bifurcation. The global dynamics in region X is characterized by three different stable limit cycles which are depicted in Fig. \ref{fig20}. These oscillatory solutions consist of two almost anti-phase and one 
 almost synchronized solutions. Two almost anti-phase solutions are recognizable throughout their different frequencies, as is depicted in Fig. \ref{fig21}. Therefore, the system as a whole is multi-stable, two non-trivial rest points and three limit cycles. As the parameter $\tau$ changes from region X to XI, the limit cycle which had become stable in region III and is of lower frequency, changes stability and becomes unstable of index 2. Hence, in this region the dynamics
could correspond to either synchronized and anti-phase spiking, or quiescence.
We should remind that, for narrow regions of delay parameter, by the mechanism stated above and in the previous section, the parameter $\tau$ can either suppress periodic spiking, or induce new ones with different frequencies and different synchrony.
\section{conclusion}\label{sec:conclusion}
As a typical representation of excitable systems, FHN
model is widely used to research some nonlinear phenomena in the nervous system.
In this paper, possible bifurcations were described in two synaptically coupled FHN neurons with the time delay. We showed that with the coupling strength changing, two coupled nonidentical neurons can exhibit rich bifurcation behaviors. We also showed that the dynamics in two coupled FHN neurons is drastically
changed due to the effect of time delay. We used theory of dynamical systems to describe typical behaviors, such as  quiescence, almost synchronized and anti-phase periodic solutions,  or coexistence of these states. It illustrates that the time delay can either suppress periodic activity or induce new ones, depending on the value of the time delay and kind of bifurcation. All possible Hopf, double Hopf,  pitchfork, pitchfork of limit cycles, homoclinic, fold of limit cycles, and torus bifurcations of our system are derived, and total bifurcation diagrams for trivial and non-trivial rest points are reported. Actually this kind of bifurcation study is a serious task for the prediction and detection of all possible dynamical behaviors. Missing some kind of solutions is possible, if we only aid simulation study. We should emphasis that our investigation provides a detailed scenario for the Hopf bifurcation diagrams appeared in similar studies.  
These results are interesting from the point of view of applications, since our generic model is representative for a
wide range of real-world systems.
In the future research we will classify possible activities for two coupled nonidentical FHN neurons without inherent symmetry.




\bibliographystyle{model1-num-names}
\bibliography{ref}







\end{document}